\input amstex.tex
\documentstyle{amsppt}
\magnification=\magstep1
\vsize 8.5truein
\hsize 6truein

\define\flow{\left(\bold{M},\{S^t\}_{t\in\Bbb R},\mu\right)}

\define\traj{S^{[a,b]}x_0}
\define\endrem{}
 
\noindent
February 5, 2003

\hfill{\it Dedicated to Ya. G. Sinai}

\hfill{\it honoring his 65th birthday}

\bigskip \bigskip

\heading
Proof of the Boltzmann-Sinai Ergodic Hypothesis \\
for Typical Hard Disk Systems
\endheading
 
\bigskip \bigskip
 
\centerline{{\bf N\'andor Sim\'anyi}
\footnote{Research supported by the National Science Foundation, grant
DMS-0098773.}}

\bigskip \bigskip

\centerline{University of Alabama at Birmingham}
\centerline{Department of Mathematics}
\centerline{Campbell Hall, Birmingham, AL 35294 U.S.A.}
\centerline{E-mail: simanyi\@math.uab.edu}

\bigskip \bigskip

\hbox{\centerline{\vbox{\hsize 8cm {\bf Abstract.} We consider the system of
$N$ ($\ge2$) hard disks of masses $m_1,\dots,m_N$ and radius $r$ in the flat
unit torus $\Bbb T^2$. We prove the ergodicity (actually, the B-mixing
property) of such systems for almost every selection $(m_1,\dots,m_N;r)$ of 
the outer geometric parameters.}}}

\bigskip \bigskip

\noindent
Primary subject classification: 37D50

\medskip

\noindent
Secondary subject classification: 34D05

\bigskip \bigskip

\heading
\S1. Introduction
\endheading

\bigskip \bigskip

Hard ball systems or, a bit more generally, mathematical billiards
constitute an important and quite interesting family of dynamical systems
being intensively studied by dynamicists and researchers of mathematical
physics, as well. These dynamical systems pose many challenging mathematical
questions, most of them concerning the ergodic (mixing) properties of such
systems. The introduction of hard ball systems and the first major steps in
their investigations date back to the 40's and 60's, see Krylov's paper
[K(1979)] and Sinai's ground-breaking works [Sin(1963)] and [Sin(1970)], 
in which the author --- among other things --- formulated the modern version
of Boltzmann's ergodic hypothesis (what we call today the Boltzmann--Sinai
ergodic hypothesis) by claiming that every hard ball system in a flat torus
is ergodic, of course after fixing the values of the trivial flow-invariant
quantities. In the papers
[Sin(1970)] and [B-S(1973)] Bunimovich and Sinai proved
this hypothesis for two hard disks on the two-dimensional unit torus
$\Bbb T^2$. The generalization of this result to higher dimensions $\nu>2$
took fourteen years, and was done by Chernov and Sinai in [S-Ch(1987)].
Although the model of two hard balls in $\Bbb T^\nu$ is already
rather involved technically, it is still a so called strictly dispersive
billiard system, i. e. such that the smooth components of the boundary
$\partial\bold Q$ of the configuration space are strictly concave from
outside $\bold Q$. (They are bending away from $\bold Q$.)
The billiard systems of more than two hard balls in
$\Bbb T^\nu$ are no longer strictly dispersive, but just semi-dispersive
(strict concavity of the smooth components of $\partial\bold Q$
is lost, merely concavity persists), and this circumstance causes a lot
of additional technical troubles in their study. In the series of my joint
papers with A. Kr\'amli and D. Sz\'asz [K-S-Sz(1989)], [K-S-Sz(1990)],
[K-S-Sz(1991)], and  [K-S-Sz(1992)]
we developed several new methods, and proved the ergodicity of 
more and more complicated semi-dispersive billiards culminating in the proof
of ergodicity of four billiard balls in the torus $\Bbb T^\nu$ 
($\nu\ge 3$), [K-S-Sz(1992)]. Then, in 1992, Bunimovich, Liverani,
Pellegrinotti and Sukhov [B-L-P-S(1992)] were able to prove the ergodicity for
some systems with an arbitrarily large number of hard balls. The shortcoming
of their model, however, is that, on one hand, they restrict the types of all
feasible ball-to-ball collisions, on the other hand, they introduce some
additional
scattering effect with the collisions at the strictly concave wall of the
container. The only result with an arbitrarily large number of balls in a
flat unit torus $\Bbb T^\nu$ was achieved in [Sim(1992-A-B)], 
where the author managed to
prove the ergodicity (actually, the K-mixing property) of $N$ hard balls in
$\Bbb T^\nu$, provided that $N\le\nu$. The annoying shortcoming of that result
is that the larger the number of balls $N$ is, larger and larger dimension
$\nu$ of the ambient container is required by the method of the proof.

On the other hand, if someone considers a hard ball system in an
elongated torus which is long in one direction but narrow in the others,
so that the balls must keep their cyclic order in the ``long direction''
(Sinai's ``pen-case'' model), then the technical difficulties can be handled,
thanks to the fact that the collisions of balls are now restricted to
neighboring pairs (in the cyclic order). The hyperbolicity of such models in 
three dimensions and the ergodicity in dimension four have been proved in 
[S-Sz(1995)]. 

The positivity of the metric entropy for several systems of hard balls
can be proved relatively easily, as was shown in the paper [W(1988)]. 
The papers [L-W(1995)] and [W(1990)] are nice surveys describing a 
general setup leading to the technical problems treated in a series of
research papers. For a comprehensive survey of the results and open problems 
in this field, see [Sz(1996)].

Pesin's theory [P(1977)] on the ergodic properties of non-uniformly hyperbolic,
smooth dynamical systems has been generalized substantially to dynamical 
systems with singularities (and with a relatively mild behavior near the
singularities) by A. Katok and J-M. Strelcyn [K-S(1986)]. Since then, the
so called Pesin's and Katok-Strelcyn's theories have become part of the 
folklore in the theory of dynamical systems. They claim that --- under some 
mild regularity conditions, particularly near the singularities --- every
non-uniformly hyperbolic and ergodic flow enjoys the Kolmogorov-mixing
property, shortly the K-mixing property.

Later on it was discovered and proved in [C-H(1996)] and [O-W(1998)] that the
above mentioned fully hyperbolic and ergodic flows with singularities turn out
to be automatically having the Bernoulli mixing (B-mixing) property. It is
worth noting here that almost every semi-dispersive billiard system, 
especially every hard ball system, enjoys those mild regularity conditions
imposed on the systems (as axioms) by [K-S(1986)], [C-H(1996)], and
[O-W(1998)]. In other words, for a hard ball flow 
$\left(\bold M,\{S^t\},\mu\right)$ the (global) ergodicity of the systems
actually implies its full hyperbolicity and the B-mixing property, as well.

Finally, in our joint venture with D. Sz\'asz [S-Sz(1999)], we prevailed
over the difficulty caused by the low value of the dimension $\nu$ by
developing a brand new algebraic approach for the study of hard ball systems.
That result, however, only establishes complete hyperbolicity (nonzero Lyapunov
exponents almost everywhere) for $N$ balls in $\Bbb T^\nu$. The ergodicity
appeared to be a harder task.

\medskip

Consider the $\nu$-dimensional ($\nu\ge2$), standard, flat, unit torus
$\Bbb T^\nu=\Bbb R^\nu/\Bbb Z^\nu$ as the vessel containing $N$ ($\ge2$)
hard balls (spheres) $B_1,\dots,B_N$ with positive masses 
$m_1,\dots,m_N$ and (just for simplicity) common radius $r>0$. We always
assume that the radius $r>0$ is not too big, so
that even the interior of the arising
configuration space $\bold Q$ is connected. Denote the center of the ball
$B_i$ by $q_i\in\Bbb T^\nu$, and let $v_i=\dot q_i$ be the velocity of the
$i$-th particle. We investigate the uniform motion of the balls
$B_1,\dots,B_N$ inside the container $\Bbb T^\nu$ with half a unit of total 
kinetic energy: $E=\dfrac{1}{2}\sum_{i=1}^N m_i||v_i||^2=\dfrac{1}{2}$.
We assume that the collisions between balls are perfectly elastic. Since
--- beside the kinetic energy $E$ --- the total momentum
$I=\sum_{i=1}^N m_iv_i\in\Bbb R^\nu$ is also a trivial first integral of the
motion, we make the standard reduction $I=0$. Due to the apparent translation
invariance of the arising dynamical system, we factorize out the configuration
space with respect to uniform spatial translations as follows:
$(q_1,\dots,q_N)\sim(q_1+a,\dots,q_N+a)$ for all translation vectors
$a\in\Bbb T^\nu$. The configuration space $\bold Q$ of the arising flow
is then the factor torus 
$\left(\left(\Bbb T^\nu\right)^N/\sim\right)\cong\Bbb T^{\nu(N-1)}$
minus the cylinders
$$
C_{i,j}=\left\{(q_1,\dots,q_N)\in\Bbb T^{\nu(N-1)}\colon\;
\text{dist}(q_i,q_j)<2r \right\}
$$
($1\le i<j\le N$) corresponding to the forbidden overlap between the $i$-th
and $j$-th spheres. Then it is easy to see that the compound configuration
point
$$
q=(q_1,\dots,q_N)\in\bold Q=\Bbb T^{\nu(N-1)}\setminus
\bigcup_{1\le i<j\le N}C_{i,j}
$$
moves in $\bold Q$ uniformly with unit speed and bounces back from the
boundaries $\partial C_{i,j}$ of the cylinders $C_{i,j}$ according to the
classical law of geometric optics: the angle of reflection equals the angle of
incidence. More precisely: the post-collision velocity $v^+$ can be obtained
from the pre-collision velocity $v^-$ by the orthogonal reflection across the
tangent hyperplane of the boundary $\partial\bold Q$ at the point of collision.
Here we must emphasize that the phrase ``orthogonal'' should be understood 
with respect to the natural Riemannian metric (the so called mass metric)
$||dq||^2=\sum_{i=1}^N m_i||dq_i||^2$ in the configuration space $\bold Q$.
For the normalized Liouville measure $\mu$ of the arising flow
$\{S^t\}$ we obviously have $d\mu=\text{const}\cdot dq\cdot dv$, where
$dq$ is the Riemannian volume in $\bold Q$ induced by the above metric and
$dv$ is the surface measure (determined by the restriction of the
Riemannian metric above) on the sphere of compound velocities
$$
\Bbb S^{\nu(N-1)-1}=\left\{(v_1,\dots,v_N)\in\left(\Bbb R^\nu\right)^N\colon\;
\sum_{i=1}^N m_iv_i=0 \text{ and } \sum_{i=1}^N m_i||v_i||^2=1 \right\}.
$$
The phase space $\bold M$ of the flow $\{S^t\}$ is the unit tangent bundle
$\bold Q\times\Bbb S^{d-1}$ of the configuration space $\bold Q$. (We will 
always use the shorthand notation $d=\nu(N-1)$ for the dimension of the 
billiard table $\bold Q$.) We must, however, note here that at the boundary
$\partial\bold Q$ of $\bold Q$ one has to glue together the pre-collision and
post-collision velocities in order to form the phase space $\bold M$, so
$\bold M$ is equal to the unit tangent bundle $\bold Q\times\Bbb S^{d-1}$
modulo this identification.

A bit more detailed definition of hard ball systems with arbitrary masses,
as well as their role in the family of cylindric billiards, can be found in
\S4 of [S-Sz(2000)] and in \S1 of [S-Sz(1999)]. We denote the
arising flow by $\flow$.

\medskip

The joint work of Ya. G. Sinai and N. I. Chernov [S-Ch(1987)] paved the way for
further fundamental results concerning the ergodicity of $\flow$. They proved 
there a strong result on local ergodicity: An open neighborhood 
$U\subset\bold M$ of every phase point with a hyperbolic trajectory (and with
at most one singularity on its trajectory) belongs to a single ergodic 
component of the billiard flow $\flow$, of course, modulo the zero measure
sets. An immediate consequence of this result is the (hyperbolic) ergodicity
of the hard ball systems with $N=2$ and $\nu\ge2$.

\medskip

\subheading{Remark} It is worth noting here that the proof of the above 
mentioned Theorem on Local Ergodicity by Chernov and Sinai necessitates the
assumption of an annoying technical condition, the so called ``Chernov-Sinai
Ansatz'', see Condition 3.1 in [K-S-Sz(1990)]. The first part of \S3 of
this paper will be devoted for proving this condition.

\medskip

In the series of papers [K-S-Sz(1989)], [K-S-Sz(1991)], 
[K-S-Sz(1992)], [Sim(1992 -A)], and [Sim(1992-B)],
the authors developed a powerful, three-step strategy for 
proving the (hyperbolic) ergodicity of hard ball systems. First of all,
all these proofs are inductions on the number $N$ of balls involved in the
problem. Secondly, the induction step itself consists of the following three
major steps:

\medskip

\subheading{Step I} To prove that every non-singular (i. e. smooth)
trajectory segment $\traj$ with a ``combinatorially rich'' (in a well
defined sense) symbolic collision sequence is automatically sufficient
(or, in other words, ``geometrically hyperbolic'', see below in \S2), 
provided that the phase point $x_0$ does not belong to a countable union $J$
of smooth sub-manifolds with codimension at least two. (Containing the 
exceptional phase points.)

The exceptional set $J$ featuring this result is negligible in our dynamical
considerations --- it is a so called slim set. For the basic properties of
slim sets, please see \S2.7 below.

\medskip

\subheading{Step II} Assume the induction hypothesis, i. e. that all hard
ball systems with $n$ balls ($n<N$) are (hyperbolic and) ergodic. Prove that
there exists a slim set $S\subset\bold M$ 
(see \S2.7) with the following property:
For every phase point $x_0\in\bold M\setminus S$ the entire trajectory
$S^{\Bbb R}x_0$ contains at most one singularity and its symbolic collision
sequence is combinatorially rich, just as required by the result of Step I.

\medskip

\subheading{Step III} By using again the induction hypothesis, prove that
almost every singular trajectory is sufficient in the time interval
$(t_0,+\infty)$, where $t_0$ is the time moment of the singular reflection.
(Here the phrase ``almost every'' refers to the volume defined by the induced
Riemannian metric on the singularity manifolds.)

We note here that the almost sure sufficiency of the singular trajectories
(featuring Step III) is an essential condition for the proof of the celebrated
Theorem on Local Ergodicity for algebraic semi-dispersive billiards
proved by B\'alint--Chernov--Sz\'asz--T\'oth in [B-Ch-Sz-T (2002)].
Under this assumption that theorem states that in any algebraic
semi-dispersive billiard system (i. e. in a system such that the smooth
components of the boundary $\partial\bold Q$ are algebraic hypersurfaces)
a suitable, open neighborhood $U_0$ of any sufficient phase point
$x_0\in\bold M$ (with at most one singularity on its trajectory) belongs
to a single ergodic component of the billiard flow $\flow$.

In an inductive proof of ergodicity, steps I and II together ensure that
there exists an arc-wise connected set
$C\subset\bold M$ with full measure, such that every phase point $x_0\in C$
is sufficient with at most one singularity on its trajectory. Then the cited
Theorem on Local Ergodicity
states that for every phase point $x_0\in C$ an open neighborhood $U_0$ of
$x_0$ belongs to one ergodic component of the flow. Finally, the connectedness
of the set $C$ and $\mu(\bold M\setminus C)=0$ easily imply that the flow
$\flow$ (now with $N$ balls) is indeed ergodic, and actually fully hyperbolic,
as well.

In the papers [K-S-Sz(1991)], [K-S-Sz(1992)] the authors followed
the strategy outlined above and obtained the (hyperbolic) ergodicity of three
and four hard balls, respectively. Technically speaking, in those papers we
always assumed tacitly that the masses of balls are equal.

The twin papers [Sim(1992-A-B)] of mine brought new topological and geometric
tools to attack the problem of ergodicity. Namely, in [Sim(1992-A)], 
a brand new topological method was developed, and that resulted in settling 
Step II of the induction, once forever. 

In the subsequent paper [Sim(1992-B)] a new combinatorial approach for 
handling Step I was developed in the case when the dimension $\nu$ of the 
toroidal container is not less than the number of balls $N$. This proves the 
ergodicity of every hard ball system with $\nu\ge N$.

\bigskip

The main result of this paper is our

\medskip

\proclaim{Theorem} In the case $\nu=2$ (i. e. for hard disks in $\Bbb T^2$)
for almost every selection $(r;\, m_1,\dots,m_N)$ of the outer geometric
parameters from the region $0<r<r_0$, $m_i>0$, (here the 
inequality $r<r_0$ just describes the region where the interior of the 
configuration space is connected) it is true that the billiard flow
$\left(\bold M_{\vec m,r},\{S^t\},\mu_{\vec m,r}\right)$ of the $N$-disk
system is ergodic and completely hyperbolic. Then, following from the results
of Chernov--Haskell [C-H(1996)] and Ornstein--Weiss [O-W(1998)], such a
semi-dispersive billiard system actually enjoys the B-mixing property, as well.

\endproclaim

\medskip

A few remarks concerning this theorem are now in place. 

\medskip

\subheading{Remark 1} The above inequality $r<r_0$ corresponds to physically 
relevant situations. Indeed, in the case $r\ge r_0$ the particles would not 
have enough room even to freely exchange positions. 

\medskip

\subheading{Remark 2} Below we present an inductive proof following 
the above drafted three-step strategy
I--III amended in such a way that the exceptional set $J$ featuring Step I
is no longer a countable union of codimension-two (i. e. at least two) sets
but, rather, it is a countable union of proper (i. e. of codimension
at least one) sub-manifolds. This shortcoming of Step I makes it possible
(in principle, at least) that countably many open ergodic components
$C_1,\, C_2,\dots$ coexist in such a way that they are separated from each
other by codimension-one, smooth, exceptional sub-manifolds $J$ of $\bold M$
featuring Step I. The main contents of the present paper is to exclude this
possibility, and this is precisely what is going on in \S4--8 below. It is 
just this proof of the non-existence of separating manifolds $J$ that
essentially uses the dimension condition $\nu=2$.

\medskip

\subheading{Remark 3}
The last remark concerns the fact that --- at least in principle --- an
unspecified zero measure set of the outer geometric parameters
$(m_1,\dots,m_N;r)$ has to be dropped in the theorem. But why? The reason is
the same as for the dropping of the zero set in the main theorem of
[S-Sz(1999)], in which we proved that a hard ball system (in any given
dimension $\nu\ge2$) is almost surely fully hyperbolic, that is, its relevant
Lyapunov exponents are nonzero almost everywhere. In fact, in the proof of
Proposition 3.1 below (which is required for the proof of the Chernov-Sinai
Ansatz, i. e. Step III) we successfully applied the algebraic method developed
in [S-Sz(1999)]. Proposition 3.1 asserts that the intersection of the
exceptional set $J$ (featuring Step I) and the singularity set
$\Cal S\Cal R^+$ (see the two paragraphs preceding Proposition 3.1) has at
least two codimensions, that is, $J$ and $\Cal S\Cal R^+$ cannot even locally
coincide.

\bigskip

The paper is organized as follows. After putting forward the prerequisites
in \S2, in the subsequent section we carry out the inductive proof of the
ergodicity by assuming the non-existence of the separating manifolds $J$.
Then all remaining sections 4--8 are devoted to the proof of the
non-existence of $J$-manifolds. That proof contains a lot of new geometric
ideas. Finally, in \S9, two remarks conclude the article. One of them regards
the role of Proposition 3.1 in the entire proof, while the other one applies
the method of [S-W(1989)] to prove the striking fact that a typical
(i. e. an ergodic, or B-mixing) hard disk system retains its B-mixing property
even if one omits the translation factorization
$(q_1,\dots,q_N)\sim(q_1+a,\dots,q_N+a)$ of the configuration space, despite
the fact that the dropping of this factorization introduces 2 zero Lyapunov
exponents!

\bigskip \bigskip
 
\heading
\S 2. Prerequisites
\endheading
 
\bigskip \bigskip  

\subheading{2.1 Cylindric billiards} Consider the $d$-dimensional
($d\ge2$) flat torus $\Bbb T^d=\Bbb R^d/\Cal L$ supplied with the
usual Riemannian inner product $\langle\, .\, ,\, .\, \rangle$ inherited
from the standard inner product of the universal covering space $\Bbb R^d$.
Here $\Cal L\subset\Bbb R^d$ is assumed to be a lattice, i. e. a discrete
subgroup of the additive group $\Bbb R^d$ with $\text{rank}(\Cal L)=d$.
The reason why we want to allow general lattices, other than just the
integer lattice $\Bbb Z^d$, is that otherwise the hard ball systems would
not be covered. The geometry of the structure lattice $\Cal L$ in the
case of a hard ball system is significantly different from the geometry
of the standard lattice $\Bbb Z^d$ in the standard Euclidean space
$\Bbb R^d$, see later in this section.

The configuration space of a cylindric billiard is
$\bold Q=\Bbb T^d\setminus\left(C_1\cup\dots\cup C_k\right)$, where the
cylindric scatterers $C_i$ ($i=1,\dots,k$) are defined as follows.

Let $A_i\subset\Bbb R^d$ be a so called lattice subspace of $\Bbb R^d$,
which means that $\text{rank}(A_i\cap\Cal L)=\text{dim}A_i$. In this case
the factor $A_i/(A_i\cap\Cal L)$ is a sub-torus in $\Bbb T^d=\Bbb R^d/\Cal L$
which will be taken as the generator of the cylinder 
$C_i\subset\Bbb T^d$, $i=1,\dots,k$. Denote by $L_i=A_i^\perp$ the
ortho-complement of $A_i$ in $\Bbb R^d$. Throughout this paper we will
always assume that $\text{dim}L_i\ge2$. Let, furthermore, the numbers
$r_i>0$ (the radii of the spherical cylinders $C_i$) and some translation
vectors $t_i\in\Bbb T^d=\Bbb R^d/\Cal L$ be given. The translation
vectors $t_i$ play a crucial role in positioning the cylinders $C_i$
in the ambient torus $\Bbb T^d$. Set
$$
C_i=\left\{x\in\Bbb T^d\colon\; \text{dist}\left(x-t_i,A_i/(A_i\cap\Cal L)
\right)<r_i \right\}.
$$
In order to avoid further unnecessary complications, we always assume that
the interior of the configuration space 
$\bold Q=\Bbb T^d\setminus\left(C_1\cup\dots\cup C_k\right)$ is connected.
The phase space $\bold M$ of our cylindric billiard flow will be the
unit tangent bundle of $\bold Q$ (modulo some natural gluing at its
boundary), i. e. $\bold M=\bold Q\times\Bbb S^{d-1}$. (Here $\Bbb S^{d-1}$
denotes the unit sphere of $\Bbb R^d$.)

The dynamical system $\flow$, where $S^t$ ($t\in\Bbb R$) is the dynamics 
defined by the uniform motion inside the domain $\bold Q$ and specular
reflections at its boundary (at the scatterers), and $\mu$ is the
Liouville measure, is called a cylindric billiard flow we want to
investigate. 

We note that the cylindric billiards --- defined above --- belong to the wider
class of so called semi-dispersive billiards, which means that the smooth 
components $\partial\bold Q_i$ of the boundary $\partial\bold Q$
of the configuration space $\bold Q$ are (not necessarily strictly) concave
from outside of $\bold Q$,
i. e. they are bending away from the interior of $\bold Q$. As to the notions 
and notations in connection with semi-dispersive billiards, the reader is 
kindly referred to the paper [K-S-Sz(1990)].

\bigskip

\subheading{2.2 Hard ball systems} Hard ball systems in the standard
unit torus $\Bbb T^\nu=\Bbb R^\nu/\Bbb Z^\nu$ ($\nu\ge2$) with positive masses
$m_1,\dots,m_N$ are described (for example) in \S 1 of [S-Sz(1999)].
These are the dynamical systems describing the motion of $N$ ($\ge2$) hard
balls with radii $r_1,r_2,\dots,r_N$ and positive masses $m_1,\dots,m_N$ in 
the standard unit torus $\Bbb T^\nu=\Bbb R^\nu/\Bbb Z^\nu$. (For simplicity
we assume that these radii have the common value $r$.) The center of the
$i$-th ball is denoted by $q_i$ ($\in\Bbb T^\nu$), its time derivative is
$v_i=\dot q_i$, $i=1,\dots,N$. One uses the standard reduction of kinetic
energy $E=\frac{1}{2}\sum_{i=1}^N m_i||v_i||^2=\frac{1}{2}$.
The arising configuration space (still without the removal of the scattering
cylinders $C_{i,j}$) is the torus
$$
\Bbb T^{\nu N}=\left(\Bbb T^{\nu}\right)^N=\left\{(q_1,\dots,q_N)\colon\;
q_i\in\Bbb T^\nu,\; i=1,\dots,N\right\}
$$
supplied with the Riemannian inner product (the so called mass metric)
$$
\langle v,v'\rangle=\sum_{i=1}^N m_i\langle v_i,v'_i \rangle
\tag 2.2.1
$$
in its common tangent space $\Bbb R^{\nu N}=\left(\Bbb R^{\nu}\right)^N$.
Now the Euclidean space $\Bbb R^{\nu N}$ with the inner product (2.2.1)
plays the role of $\Bbb R^d$ in the original definition of cylindric
billiards, see \S 2.1 above.

The generator subspace $A_{i,j}\subset \Bbb R^{\nu N}$ ($1\le i<j\le N$)
of the cylinder $C_{i,j}$ (describing the collisions between the $i$-th and
$j$-th balls) is given by the equation
$$
A_{i,j}=\left\{(q_1,\dots,q_N)\in\left(\Bbb R^\nu\right)^N\colon\;
q_i=q_j \right\},
\tag 2.2.2
$$
see (4.3) in [S-Sz(2000)]. Its ortho-complement 
$L_{i,j}\subset\Bbb R^{\nu N}$ is then defined by the equation

$$
L_{i,j}=\left\{(v_1,\dots,v_N)\in\left(\Bbb R^\nu\right)^N\colon\;
v_k=0 \text{ for } k\ne i,j, \text{ and } m_iv_i+m_jv_j=0 \right\},
\tag 2.2.3
$$
see (4.4) in [S-Sz(2000)].
Easy calculation shows that the cylinder $C_{i,j}$ 
(describing the overlap of the $i$-th and $j$-th balls)
is indeed spherical and the radius of its base sphere is equal to
$r_{i,j}=2r\sqrt{\frac{m_im_j}{m_i+m_j}}$, see \S 4, especially formula
(4.6) in [S-Sz(2000)].

The structure lattice $\Cal L\subset\Bbb R^{\nu N}$ is clearly the integer
lattice $\Cal L=\Bbb Z^{\nu N}$. 

Due to the presence of an extra invariant quantity
$I=\sum_{i=1}^N m_iv_i$, one usually makes the reduction
$\sum_{i=1}^N m_iv_i=0$ and, correspondingly, factorizes the configuration
space with respect to uniform spatial translations:
$$
(q_1,\dots,q_N)\sim(q_1+a,\dots,q_N+a), \quad a\in\Bbb T^\nu.
$$
The natural, common tangent space of this reduced configuration space is then
$$
\Cal Z=\left\{(v_1,\dots,v_N)\in\left(\Bbb R^\nu\right)^N\colon\;
\sum_{i=1}^N m_iv_i=0\right\}=\left(\bigcap_{i<j}A_{i,j}
\right)^\perp=\left(\Cal A\right)^\perp
\tag 2.2.4
$$
supplied again with the inner product (2.2.1), see also (4.1) and
(4.2) in [S-Sz(2000)]. The base spaces $L_{i,j}$
of (2.2.3) are obviously subspaces of $\Cal Z$, and we take
$\tilde A_{i,j}=A_{i,j}\cap\Cal Z=P_{\Cal Z}(A_{i,j})$ as the ortho-complement
of $L_{i,j}$ in $\Cal Z$. (Here $P_{\Cal Z}$ denotes the orthogonal projection
onto the space $\Cal Z$.)

Note that the configuration space of the reduced system (with 
$\sum_{i=1}^N m_iv_i=0$) is naturally the torus
$\Bbb R^{\nu N}/(\Cal A+\Bbb Z^{\nu N})=\Cal Z/P_{\Cal Z}(\Bbb Z^{\nu N})$.

\bigskip

\subheading{2.3 Singularities and Trajectory Branches} The billiard flow
$\flow$ has two types of singularities: The first type is the so called 
tangency (or tangential reflection/collision) that takes place at a phase point
$(q,v)\in\partial\bold M$, where the velocity vector happens to lie inside
the tangent space $\Cal T_q\partial\bold Q$ of the boundary $\partial\bold Q$
of the configuration space. If a trajectory hits this type of singularity,
it still has a unique continuation, the flow is still continuous at such a
phase point, but it ceases to be differentiable. The first return map $T$
to the boundary $\partial\bold M$ is no longer even continuous in any 
neighborhood of a phase point with tangential singularity.

The second type of singularity is the case of a so called ``double collision'',
when two collisions $(i,j)$ and $(j,k)$ (sharing the same disk labeled by
$j$ here) happen to take place exactly at the same time $t_0$. 
(Typically there are no more collisions taking place at time $t_0$.)
We are going to briefly describe the discontinuity of the flow
$\{S^t\}$ caused by a double collision at time $t_0$.
Assume first that the pre--collision velocities of the particles are given.
What can we say about the possible post--collision velocities? Let us perturb
the pre--collision phase point (at time $t_0-0$) infinitesimally, so that the
collisions at $\sim t_0$ occur at infinitesimally different moments. By
applying the collision laws to the arising finite sequence of collisions, 
(the finiteness follows from Theorem 1 of [B-F-K(1998)]) we
see that the post-collision velocities are fully determined by the time--
ordering of the considered collisions. Therefore, the collection of all
possible time-orderings of these collisions gives rise to a finite family of
continuations of the trajectory beyond $t_0$. They are called the
{\bf trajectory branches}. It is quite clear that similar statements can be
said regarding the evolution of a trajectory through a multiple collision
{\bf in reverse time}. Furthermore, it is also obvious that for any given
phase point $x_0\in\bold M$ there are two, $\omega$-high trees
$\Cal T_+$ and $\Cal T_-$ such that $\Cal T_+$ ($\Cal T_-$) describes all the
possible continuations of the positive (negative) trajectory
$S^{[0,\infty)}x_0$ ($S^{(-\infty,0]}x_0$). (For the definitions of trees and
for some of their applications to billiards, cf. the beginning of \S 5
in [K-S-Sz(1992)].) It is also clear that all possible continuations
(branches) of the whole trajectory $S^{(-\infty,\infty)}x_0$ can be uniquely
described by all possible pairs $(B_-,B_+)$ of infinite branches of
the trees $\Cal T_-$ and $\Cal T_+$ ($B_-\subset\Cal T_-, B_+\subset
\Cal T_+$).

Since, in the case of double collisions, there is no unique continuation of
the trajectories, we need to make a clear distinction between the set of 
reflections $\Cal S\Cal R^+$ supplied with the outgoing velocity $v^+$, and 
the set of reflections $\Cal S\Cal R^-$ supplied with the incoming velocity 
$v^-$. For typical phase points $x^+\in\Cal S\Cal R^+$ the forward trajectory 
$S^{[0,\infty)}x^+$ is non-singular and uniquely defined, and analogous 
statement holds true for typical phase points $x^-\in\Cal S\Cal R^-$ and the 
backward trajectory $S^{(-\infty,0]}x^-$. For a more detailed exposition of
singularities, the reader is kindly referred to \S2 of [K-S-Sz(1990)].

Finally, we note that the trajectory of the phase point $x_0$ has exactly two
branches, provided that $S^tx_0$ hits a singularity for a single value
$t=t_0$, and the phase point $S^{t_0}x_0$ does not lie on the intersection
of more than one singularity manifolds. (In this case we say that the 
trajectory of $x_0$ has a ``single singularity''.)

\bigskip

\subheading{2.4 Neutral Subspaces, Advance, and Sufficiency}
Consider a {\bf nonsingular} trajectory segment $S^{[a,b]}x$.
Suppose that $a$ and $b$ are {\bf not moments of collision}.

\medskip

\proclaim{Definition 2.4.1} The neutral space $\Cal N_0(S^{[a,b]}x)$
of the trajectory segment $S^{[a,b]}x$ at time zero ($a<0<b$)\ is
defined by the following formula:
$$
\aligned
\Cal N_0(S^{[a,b]}x)=\big \{W\in\Cal Z\colon\;\exists (\delta>0) \;
\text{s. t.} \; \forall \alpha \in (-\delta,\delta) \\
p\left(S^a\left(Q(x)+\alpha W,V(x)\right)\right)=p(S^ax)\text{ and }
p\left(S^b\left(Q(x)+\alpha W,V(x)\right)\right)=p(S^bx)\big\},
\endaligned
$$
where $p(Q,V)=:V$ is the projection onto the velocity sphere for any 
$(Q,V)\in\bold M$.
\endproclaim 
\endrem

(The formula for the tangent space $\Cal Z$ can be found in (2.2.4).)

It is known (see (3) in \S 3 of [S-Ch (1987)]) that
$\Cal N_0(S^{[a,b]}x)$ is a linear subspace of $\Cal Z$ indeed, and
$V(x)\in \Cal N_0(S^{[a,b]}x)$. The neutral space $\Cal N_t(S^{[a,b]}x)$
of the segment $S^{[a,b]}x$ at time $t\in [a,b]$ is defined as follows:
$$
\Cal N_t(S^{[a,b]}x)=\Cal N_0\left(S^{[a-t,b-t]}(S^tx)\right).
\tag 2.4.2
$$
It is clear that the neutral space $\Cal N_t(S^{[a,b]}x)$ can be
canonically
identified with $\Cal N_0(S^{[a,b]}x)$ by the usual identification of the
tangent spaces of $\bold Q$ along the trajectory $S^{(-\infty,\infty)}x$
(see, for instance, \S 2 of [K-S-Sz(1990)]).

Finally, the neutral space $\Cal N_0\left(S^{[a,\infty)}x\right)$ of an
unbounded trajectory segment $S^{[a,\infty)}x$ is defined as the limiting
space 

$$
\lim_{b\to\infty}\Cal N_0\left(S^{[a,b]}x\right)=
\bigcap\left\{\Cal N_0\left(S^{[a,b]}x\right)\big|\; b>a\right\},
$$
and the definitions of $\Cal N_0\left(S^{(-\infty,b]}x\right)$ and
$\Cal N_0\left(S^{(-\infty,\infty)}x\right)$ are analogous limits.

Our next  definition is  that of  the {\bf advance}. Consider a
non-singular orbit segment $S^{[a,b]}x$ with the symbolic collision sequence
$\Sigma=(\sigma_1, \dots, \sigma_n)$ ($n\ge 1$). This means the following:
For $k=1,\dots,n$ the symbol $\sigma_k=\left\{i_k,\,j_k\right\}$ 
($1\le i_k<j_k\le N$) is an unordered pair of disk labels, so that all the
collisions on the trajectory segment take place (in time ordering) between the
two disks listed in $\sigma_1,\dots,\sigma_n$, respectively. We also use the
notation $t_k$ ($a<t_1<\dots<t_n<b$) for the time moment of the $k$th collision
$\sigma_k$ on $S^{[a,b]}x$. For $x=(Q,V)\in\bold M$ and $W\in\Cal Z$, 
$\Vert W\Vert$ sufficiently small, denote $T_W(Q,V):=(Q+W,V)$.

\proclaim{Definition 2.4.3}
For any $1\le k\le n$ and $t\in[a,b]$, the advance
$$
\alpha(\sigma_k)\colon\;\Cal N_t(S^{[a,b]}x) \rightarrow  \Bbb R
$$
of the collision $\sigma_k$ is the unique linear extension of the linear 
functional $\alpha(\sigma_k)$
defined in a sufficiently small neighborhood of the origin of 
$\Cal N_t(S^{[a,b]}x)$ in the following way:
$$
\alpha(\sigma_k)(W):= t_k(x)-t_k(S^{-t}T_WS^tx).
$$
\endproclaim \endrem

It is now time to bring up the basic notion of {\bf sufficiency} 
(or, sometimes it is also called hyperbolicity) of a
trajectory (segment). This is the utmost important necessary condition for
the proof of the fundamental theorem for semi-dispersive billiards, see
Condition (ii) of Theorem 3.6 and Definition 2.12 in [K-S-Sz(1990)].

\medskip

\proclaim{Definition 2.4.4}
\roster
\item
The nonsingular trajectory segment $S^{[a,b]}x$ ($a$ and $b$ are supposed not
to be moments of collision) is said to be {\bf sufficient} if and only if
the dimension of $\Cal N_t(S^{[a,b]}x)$ ($t\in [a,b]$) is minimal, i.e.
$\text{dim}\ \Cal N_t(S^{[a,b]}x)=1$.
\item
The trajectory segment $S^{[a,b]}x$ containing exactly one singularity (a so 
called ``single singularity'', see above) is said to be {\bf sufficient} if 
and only if both branches of this trajectory segment are sufficient.
\endroster
\endproclaim \endrem

\medskip

\proclaim{Definition 2.4.5}
The phase point $x\in\bold M$ with at most one singularity is said
to be sufficient if and only if its whole trajectory $S^{(-\infty,\infty)}x$
is sufficient, which means, by definition, that some of its bounded
segments $S^{[a,b]}x$ are sufficient.
\endproclaim \endrem

In the case of an orbit $S^{(-\infty,\infty)}x$ with a single
singularity, sufficiency means that both branches of
$S^{(-\infty,\infty)}x$ are sufficient.

\bigskip

\subheading{2.5 No accumulation (of collisions) in finite time} 
By the results of Vaserstein [V(1979)], Galperin [G(1981)] and
Burago-Ferleger-Kononenko [B-F-K(1998)], in a 
semi-dis\-per\-sive billiard flow there can only be finitely many 
collisions in finite time intervals, see Theorem 1 in [B-F-K(1998)]. 
Thus, the dynamics is well defined as long as the trajectory does not hit 
more than one boundary components at the same time.

\bigskip

\subheading{2.6 Collision graphs} Let $S^{[a,b]}x$ be a nonsingular, finite
trajectory segment with the collisions $\sigma_1,\dots,\sigma_n$ listed in
time order. (Each $\sigma_k$ is an unordered pair $(i,j)$ of different labels
$i,j\in\{1,2,\dots,N\}$.) The graph $\Cal G=(\Cal V,\Cal E)$ with vertex set
$\Cal V=\{1,2,\dots,N\}$ and set of edges 
$\Cal E=\{\sigma_1,\dots,\sigma_n\}$ is called the {\bf collision graph}
of the orbit segment $S^{[a,b]}x$.

\bigskip

\subheading{2.7 Slim sets} 
We are going to summarize the basic properties of codimension-two subsets $A$
of a smooth manifold $M$. Since these subsets $A$ are just those
negligible in our dynamical discussions, we shall call them {\bf slim}. 
As to a  broader exposition of the issues, see [E(1978)] or \S 2 of
[K-S-Sz(1991)].

Note that the dimension $\dim A$ of a separable metric space $A$ is one of the
three classical notions of topological dimension: the covering (\v
Cech-Lebesgue), the small inductive (Menger-Urysohn), or the large inductive 
(Brouwer-\v Cech) dimension. As it is known from general topology, 
all of them are the same for separable metric spaces.

\medskip

\proclaim{Definition 2.7.1}
A subset $A$ of $M$ is called slim if and only if $A$ can be covered by a 
countable family of codimension-two (i. e. at least two) closed sets of
$\mu$--measure zero, where $\mu$ is a smooth measure on $M$. (Cf.
Definition 2.12 of [K-S-Sz(1991)].)
\endproclaim \endrem

\medskip

\proclaim{Property 2.7.2} The  collection of all slim subsets of $M$ is a
$\sigma$-ideal, that is, countable unions of slim sets and arbitrary
subsets of slim sets are also slim.
\endproclaim

\medskip

\proclaim{Lemma 2.7.3} A subset $A\subset M$ is slim if and only if for
every $x\in A$ there exists an open neighborhood $U$ of $x$ in $M$ such that 
$U\cap A$ is slim. (Locality, cf. Lemma 2.14 of [K-S-Sz(1991)].)
\endproclaim

\medskip

\proclaim{Property 2.7.4} A closed subset $A\subset M$ is slim if and only
if $\mu(A)=0$ and $\dim A\le\dim M-2$.
\endproclaim

\medskip

\proclaim{Property 2.7.5 (Integrability)}
If $A\subset M_1\times M_2$ is a closed subset of the product of two manifolds,
and for every $x\in M_1$ the set
$$
A_x=\{ y\in M_2\colon\; (x,y)\in A\}
$$
is slim in $M_2$, then $A$ is slim in $M_1\times M_2$.
\endproclaim

\medskip

The following lemmas characterize the codimension-one and codimension-two sets.

\proclaim{Lemma 2.7.6} For any closed subset $S\subset M$ the following three
conditions are equivalent:

\roster

\item"{(i)}" $\dim S\le\dim M-2$;

\item"{(ii)}"  $\text{int}S=\emptyset$ and for every open connected set 
$G\subset M$ the difference set $G\setminus S$ is also connected;

\item"{(iii)}" $\text{int}S=\emptyset$ and for every point $x\in M$ and for any
open neighborhood $V$ of $x$ in $M$ there exists a smaller open neighborhood
$W\subset V$ of the point $x$ such that for every pair of points 
$y,z\in W\setminus S$ there is a continuous curve $\gamma$ in the set 
$V\setminus S$ connecting the points $y$ and $z$.

\endroster

\endproclaim

\noindent
(See Theorem 1.8.13 and Problem 1.8.E of [E(1978)].)

\medskip

\proclaim{Lemma 2.7.7} For any subset $S\subset M$ the condition 
$\dim S\le\dim M-1$ is equivalent to $\text{int}S=\emptyset$.
(See Theorem 1.8.10 of [E(1978)].)
\endproclaim

\medskip

We recall an elementary, but important lemma (Lemma 4.15 of [K-S-Sz(1991)]).
Let $R_2$ be the set of phase points 
$x\in\bold M\setminus\partial\bold M$ such that the trajectory 
$S^{(-\infty,\infty)}x$ has more than one singularities.

\proclaim{Lemma 2.7.8} The set $R_2$ is a countable union of codimension-two
smooth sub-manifolds of $M$ and, being such, it is slim.
\endproclaim

\medskip

The next lemma establishes the most important property of slim sets which
gives us the fundamental geometric tool to connect the open ergodic components
of billiard flows.

\proclaim{Lemma 2.7.9}
If $M$ is connected, then the complement $M\setminus A$
of a slim set $A\subset M$ necessarily contains an arc-wise connected,
$G_\delta$ set of full measure. (See Property 3 of \S 4.1 in [K-S-Sz(1989)].
The  $G_\delta$ sets are, by definition, the countable intersections
of open sets.)
\endproclaim

\bigskip \bigskip

\heading
\S 3. The Inductive Proof of the Theorem \\
(Using the results of \S 4--8)
\endheading

\bigskip \bigskip

In this section we prove our theorem by using an induction on the number of
disks $N$ ($\ge2$). Consider therefore an $N$-disk billiard flow
$$
\left(\bold M_{\vec m,r},\{S^t\}_{t\in\Bbb R},\mu_{\vec m,r}\right)=\flow
$$ 
in the standard unit $2$-torus $\Bbb T^2=\Bbb R^2/\Bbb Z^2$ with the
$N+1$-tuple of outer geometric parameters $(m_1,\dots,m_N;r)$, for which 
even the interior of the phase space is connected, see the previous section.

As Corollary 3.24 and Lemma 4.2 of [Sim(2002)] state, there exists
a positive integer $C(N)$ with the following property: If the non-singular
trajectory segment $S^{[0,T]}x_0$ of $\flow$ has at least $C(N)$ consecutive,
connected collision graphs, then there exists an open neighborhood $U_0$
of $x_0$ in $\bold M$ and a proper (i. e. of codimension at least one)
algebraic set $N_0\subset U_0$ such that $S^{[0,T]}y$ is sufficient 
(or, geometrically hyperbolic, see \S2) for all $y\in U_0\setminus N_0$.

\medskip

Consider the $(2d-3)$-dimensional, compact cell complex $\Cal S\Cal R^+$
of singular reflections $x=(q,v^+)\in\partial\bold M$ supplied with the
outgoing (post collision) velocity $v^+$, so that the positive orbit is well
defined, as long as there is no other singularity on $S^{(0,\infty)}x$.
Recall that, as it follows from Lemma 4.1 of [K-S-Sz(1990)], the set of phase
points with more than one singularities is a countable union of smooth
sub-manifolds with codimension at least two (see (2.7.8) above),
thus the positive 
semi-trajectory $S^{(0,\infty)}x$ is non-singular for $\nu$-almost every
$x\in\Cal S\Cal R^+$, where $\nu$ denotes the hypersurface measure on
$\Cal S\Cal R^+$. Also recall that $2d-1=4N-5$ is the dimension of the phase
space $\bold M$.

Let $\Sigma=(\sigma_1,\dots,\sigma_n)$ be any fixed symbolic collision
sequence with at least $C(N)$ consecutive, connected collision graphs, i. e.
a so called $C(N)$-rich symbolic sequence in the sense of Key Lemma 4.1 and
Remark 4.1/b of [S-Sz(1999)]. Let, moreover, $x_0\in\Cal S\Cal R^+$ be an
element of a $(2d-3)$-dimensional, open cell $C$ of $\Cal S\Cal R^+$
(that is, $x_0$ does not belong to the $(2d-4)$-skeleton of the
$(2d-3)$-dimensional cell complex $\Cal S\Cal R^+$) with a non-singular
trajectory segment $S^{[0,T]}x_0$, for which the symbolic collision
sequence is the given $\Sigma$. We will need the following generalization of
Key Lemma 4.1 of [S-Sz(1999)], which claims that the exceptional algebraic
set $N_0$ from above (containing all phase points $y\in U_0$ around $x_0$
for which $S^{[0,T]}y$ is not sufficient) cannot even locally coincide
with the invariant hull of $\Cal S\Cal R^+$.

\medskip

Here we briefly recall the definition of a (compact) cell complex $K_n$ 
(of dimension $n$) from topology. It is defined by an induction on the 
dimension $n$, which is a non-negative integer. A zero-dimensional cell 
complex $K_0$ is a finite space with the discrete topology. For a positive 
integer $n$, an $n$-dimensional cell complex $K_n$ is a compact, metrizable 
space with a given closed subset $K_{n-1}$ (called the ($n-1$)-skeleton of 
$K_n$) so that $K_n\setminus K_{n-1}=\bigcup_{i=1}^k C_i$, where 
$C_1,\,\dots,\,C_k$ ($k\ge1$) are mutually disjoint open sets (the so called 
$n$-cells) supplied with a continuous map $\phi_i:\; \overline{B}^n\to K_n$, 
so that

\medskip

(i) $\phi_i\big|B^n:\; B^n\to C_i$ is a homeomorphism;

\medskip

(ii) $\phi_i$ maps the boundary $\partial B^n$ into the ($n-1$)-skeleton
$K_{n-1}$ (the so called gluing map).

\medskip

\noindent
(Here $\overline{B}^n$ is the closed, unit $n$-ball, and $B^n$ is its 
interior. The map $\phi_i\big|B^n:\; B^n\to C_i$ is often called the 
coordinate chart for the open $n$-cell $C_i$.)

In our examples the coordinate charts of the cells are smooth, and our compact
cell complexes turn out to be finite unions of smooth submanifolds of a
euclidean space.

\medskip

\proclaim{Proposition 3.1} For almost every $(N+1)$-tuple
$(m_1,\dots,m_N;r)$ of the outer geometric parameters the set
$$
\left\{y\in U_0\cap C\colon\; S^{[0,T]}y \text{ is not sufficient}\right\}
\tag 3.2
$$
has an empty interior in $C$. (It is actually a finite union of proper,
real analytic subsets of $C$.)

\endproclaim

\medskip

\subheading{Proof (A brief outline)} Since this generalization of 
Key Lemma 4.1 of [S-Sz(1999)] is a direct application of the proof of that
lemma (in which application all steps of the mentioned proof need to be
repeated with only minor changes), hereby we will only briefly sketch the
proof by mainly shedding light on the important steps, during which we point
out the differences between the original proof of Key Lemma 4.1 and its 
modification that proves Proposition 3.1 above. This sketch of the proof
will be subdivided into $12$ points, as follows.

\medskip

$1^o$ In order to facilitate the use of arithmetics for the kinetic variables,
we lift the entire system to the universal covering space $\Bbb R^2$ of
$\Bbb T^2$ by introducing the notion of adjustment vectors, see Propositions
3.1 and 3.3 of [S-Sz(1999)].

\medskip

$2^o$ We need to complexify the system, to introduce the algebraically
independent initial variables, the polynomial equations defining the
algebraic dynamics, the algebraic functions in terms of the initial
kinetic variables, and the tower of fields made up by all kinetic variables
of orbit segments with a symbolic collision sequence 
$(\sigma_1,\dots,\sigma_n;\, a_1,\dots,a_n)$, along the lines of 
pp. 49--54 of [S-Sz(1999)].

\medskip

$3^o$ We should introduce the complex neutral space $\Cal N(\omega)$,
just as in (3.21) of [S-Sz(1999)]. By dropping the factorization with
respect to uniform spatial translations, the condition of sufficiency now
becomes $\text{dim}\Cal N(\omega)=\nu+1=3$.

\medskip

$4^o$ Just for technical reasons, the reciprocal $1/L$ of the size $L$ of the
container $\Bbb T_L^2=\Bbb R^2/L\cdot\Bbb Z^2$ (containing disks with unit 
radius) is replaced by the radius of disks $r$ moving in the unit torus 
$\Bbb R^2/\Bbb Z^2$.

\medskip

$5^o$ By using the defining equation of the actual singularity
cell $C$, we eliminate one variable out of the initial ones to gain
again algebraic independence, despite considering singular phase points on
$C$. We express the sufficiency of $\omega=S^{[0,T]}x$ ($x\in C$) in terms of
the remaining, algebraically independent initial variables. Non-sufficiency 
again proves to be equivalent to the simultaneous vanishing of finitely many 
polynomials, in the spirit of Lemma 4.2 from [S-Sz(1999)].

\medskip

$6^o$ One reformulates the claim of Proposition 3.1 in terms of the initial
kinetic variables. The negation of that assertion proves to be the identical
vanishing of certain algebraic functions, see also Lemma 4.2 and its proof
in [S-Sz(1999)].

\medskip

$7^o$ Use Property (A) (the technical property defined in 3.31 of [S-Sz(1999)])
and the concept of combinatorial richness of the 
symbolic sequence of $S^{[0,T]}x$ (of containing at least $C(N)$ consecutive,
connected collision graphs), just like in Key Lemma 4.1 and Remark 4.1/b 
in [S-Sz(1999].

\medskip

$8^o$ Carry out an inductive proof for Proposition 3.1 above. The induction
goes on with respect to the number of disks $N$ ($\ge2$), and this induction
is independent of the outer induction to be carried out to prove the Theorem.
The statement is obviously true for $N=2$. We assume $N>2$ and the induction
hypothesis, and perform an indirect proof for the induction step by assuming
the negation of Proposition 3.1 for the complexified $N$-disk system. By using
the combinatorial richness formulated in Key Lemma 4.1 of [S-Sz(1999)],
one selects a label $i\in\{1,2,\dots,N\}$ for the substitution $m_i=0$,
along the lines of Lemma 4.43 of [S-Sz(1999]. The substitution $m_i=0$
results in a derived scheme $(\Sigma',\Cal A')$ by also preserving
Property (A), see Definition 4.11, Main Lemma 4.21, Remark 4.22, and Corollary
4.35 of [S-Sz(1999].

\medskip

$9^o$ Describe non-sufficiency in the case $m_i=0$ along the lines of 
Lemma 4.9 from [S-Sz(1999].

\medskip

$10^o$ From the indirect assumption one obtains that the induction hypothesis
is false for the $(N-1)$-disk system, just like in Lemma 4.40 of [S-Sz(1999)].

\medskip

$11^o$ From the complex version of the analogue of Key Lemma 4.1 one switches
to the real case, just as in the fourth paragraph on page 88 of [S-Sz(1999].

\medskip

$12^o$ From the real version of the analogue of Key Lemma 4.1 one obtains
Proposition 3.1 of this article by dropping a null set of $(N+1)$-tuples
$(m_1,\dots,m_N;r)$ of outer geometric parameters, precisely as in the first 
paragraph of page 93 of [S-Sz(1999)]. \qed

\bigskip

By the results of Vaserstein [V(1979)], Galperin [G(1981)], and 
Burago-Ferleger-Kononenko [B-F-K(1998)], in a semi-dispersive billiard flow
there can only be finitely many collisions in finite time, see Theorem 1
of [B-F-K(1998)], see also 2.5 above. Thus the dynamics is well defined, as 
long as the trajectory does not hit more than one boundary components at the 
same time.

\medskip

Lemma 4.1 of [K-S-Sz(1990)] claims that the set

$$
\Delta_2=\left\{x\in\bold M\colon\;\exists\text
{ at least $2$ singularities on } S^{\Bbb R}x\right\}
$$
is a countable union of smooth sub-manifolds of $\bold M$ with codimension at
least two. Especially, the set $\Delta_2$ is slim, i. e. negligible in our
considerations, see also Lemma 2.7.8 above.

\medskip

By using the results of \S 4--8, we are now going
to prove the theorem by an induction on the number of disks $N$ ($\ge2$).
For $N=2$ the result is proved by Sinai in [Sin(1970)].

Suppose now that $N>2$, and the theorem has been proved for every number of
disks $N'<N$. Theorem 5.1 of [Sim(1992-A)] together with the slimness of 
the set $\Delta_2$ of doubly singular phase points assert that there exists
a slim subset $S_1\subset\bold M$ of the phase space such that for every
$x\in\bold M\setminus S_1$ the phase point $x$ has at most one singularity
on its trajectory $S^{\Bbb R}x$, and each branch of
$S^{\Bbb R}x$ contains infinitely many 
consecutive, connected collision graphs. By Corollary 3.24 and Lemma 4.2 of 
[Sim(2002)], there exists a locally finite (and, therefore, countable) family
of codimension-one, smooth, exceptional sub-manifolds $J_i\subset\bold M$  
such that for every phase point $x\not\in\bigcup_{i}J_i\cup S_1$
the trajectory $S^{\Bbb R}x$ has at most one singularity and it is 
sufficient. According to the celebrated Theorem on Local Ergodicity for
algebraic semi-dispersive billiards by B\'alint--Chernov--Sz\'asz--T\'oth
[B-Ch-Sz-T (2002)] (see also Theorem 5 in [S-Ch(1987)] and Corollary 3.12 in
[K-S-Sz(1990)]) an open neighborhood $U_x\ni x$ of such a phase point
$x\not\in\bigcup_{i}J_i\cup S_1$ belongs to one ergodic component of the 
billiard flow $\flow$, therefore this billiard flow has (at most countably 
many) open ergodic components $C_1,C_2,\dots$. We note that 
Theorem 5.1 of [Sim(1992-A)] uses the induction hypothesis!

We carry out an indirect proof for the induction step. Assume that, contrary
to the assertion of the theorem, the number of the (open) ergodic components
$C_1,C_2,\dots$ is at least two. The main question is how different ergodic
components $C_i$ and $C_j$ can be separated in $\bold M$?

The above argument showed that, in the case of more than one ergodic components
$C_i$, there must exist a codimension-one, smooth (actually, analytic) 
exceptional sub-manifold $J\subset\bold M\setminus\partial\bold M$, 
a non-singular phase point $x_0\in J$, and an open ball neighborhood 
$B_0=B_0(x_0,\epsilon_0)\subset\bold M\setminus\partial\bold M$ of $x_0$
with the following properties:

\medskip

(0) The pair of sets $(B_0,B_0\cap J)$ is analytically
diffeomorphic to the standard pair $(\Bbb R^{2d-1},\Bbb R^{2d-2})$, where 
$2d-1=\text{dim}\bold M$, and the two connected components $B^1$ and $B^2$ 
of $B_0\setminus J$ belong to distinct ergodic components $C_i$ and $C_j$;

(1) For every $x\in B_0$ the semi-orbit $S^{[0,\infty)}x$ is sufficient
(hyperbolic) if and only if $x\not\in J$;

(2) The dimension $\text{dim}\Cal N_0(S^{[0,\infty)}x_0)$ of the neutral space 

$$
\Cal N_0\left(S^{[0,\infty)}x_0\right)=\bigcap_T
\Cal N_0\left(S^{[0,T]}x_0\right)
$$ 
of the semi-orbit $S^{[0,\infty)}x_0$ (see \S 2) is the minimum
possible value for all separating manifolds $J$ and phase points $x\in J$.
Then, by the upper semi-continuity of this dimension, we can assume that the
neighborhood $B_0$ is already small enough to ensure that

$$
\text{dim}\Cal N_0(S^{[0,\infty)}x)=\text{dim}\Cal N_0(S^{[0,\infty)}x_0)
\text{ for all } x\in J\cap B_0;
$$

(3) There exists a countable union $S_J$ of proper (i. e. of codimension at
least one) sub-manifolds of $J$, such that for every $x\in J\setminus S_J$
the positive orbit $S^{[0,\infty)}x$ is non-singular, and 
$x_0\in J\setminus S_J$.

\medskip

We note here the simple way how this last property can be achieved. Lemma 4.1
of [K-S-Sz(1990)] claims that two singularity manifolds corresponding
to different singularities on a trajectory can not even locally coincide.
(Their intersection has codimension at least $2$.) Thus, there are two 
possibilities for our exceptional manifold $J$: Either the set $S_J$
of phase points $x\in J$ with a singular forward orbit $S^{[0,\infty)}x$ has 
an empty interior in $J$, or not. In the first case we are done, for in that
case the subset $S_J$ of $J$ is actually a countable union of some proper, 
smooth submanifolds of $J$. In the second case, however, a small open set
$\emptyset\ne G\subset J$ happens to have the property that every phase point
$x\in G$ experiences a (simple) singularity at time $t(x)>0$ on its forward
orbit, the time moment $t(x)$ being a smooth function of $x$. Then, with some
value $t_0>\sup\left\{t(x)|\; x\in G\right\}$, we can take $S^{t_0}(G)=J_0$ 
as a new exceptional manifold in such a way that the mapping $S^{t_0}$ be
smooth on $G$ by taking the appropriate trajectory branches of $S^{t_0}x$
for $x\in G$. After switching to $J_0$ from $J$, almost every phase point
$x\in J_0$ will have no singularity on its forward orbit, according to
Lemma 4.1 of [K-S-Sz(1990)].

\medskip

We say that $J$ is a ``separating manifold''. The results of \S 4--8
assert that such a separating manifold $J$ does not exist. This contradiction
finishes the proof of the theorem.

\bigskip \bigskip

\heading
\S 4. Non-Existence of Separating $J$-Manifolds. \\
Part A: The Neutral Sector Opens Up
\endheading

\bigskip \bigskip

As we have seen in \S3, the only obstacle on the road of successfully
proving (by induction) the ergodicity of almost every hard disk system
$\flow$ is the following situation: There exists a codimension-one, smooth
(actually, analytic) sub-manifold $J\subset\bold M\setminus\partial\bold M$, 
a phase point $x_0\in J$, and an open ball neighborhood 
$B_0=B_0(x_0,\epsilon_0)\subset\bold M\setminus\partial\bold M$ of $x_0$
in $\bold M$ with the properties (0)---(3) listed at the end of the previous
section.

The assumed minimality of 

$$
\text{dim}\Cal N_0(S^{[0,\infty)}x)=\text{dim}\Cal N_0(S^{[0,\infty)}x_0)
\text{ for all } x\in J\cap B_0
\tag 4.1
$$
will have profound geometric consequences in the upcoming sections.

\medskip

First of all, we need to introduce a few notions and notations. Let
$w_0\in\Cal N_0(S^{[0,\infty)}x_0)$ be a unit neutral vector of $x_0$ with
the additional property $\langle w_0,v_0\rangle=0$, where $v_0$ is the 
velocity component of the phase point $x_0=(q_0,v_0)$. For any pair of real
numbers $(\tau_1,\tau_2)$ ($|\tau_i|<\epsilon_1$, $\epsilon_1>0$ is fixed,
chosen sufficiently small) we define $T_{\tau_1,\tau_2}x_0$ as the phase point

$$
T_{\tau_1,\tau_2}x_0=\left(q_0+\tau_1w_0,\, (1+\tau_2^2)^{-1/2}
(v_0+\tau_2w_0)\right).
\tag 4.2
$$
It follows immediately from the properties of the exceptional manifold $J$
that

$$
T_{\tau_1,\tau_2}x_0\in J\cap B_0 \text{ for } |\tau_i|<\epsilon_1,
\tag 4.3
$$
as long as the upper bound $\epsilon_1$ is selected small enough.

\bigskip

\heading
Basic Properties of $S^{[0,\infty)}x$ for $x=T_{\tau_1,\tau_2}x_0$
\endheading

\medskip

We want to investigate the positive semi-trajectories $S^{[0,\infty)}x$
of the phase points $x=T_{\tau_1,\tau_2}x_0\in J\cap B_0$,
$|\tau_i|<\epsilon_1$. The key point in their investigation is that both
translations by the vectors
$\tau_1w_0$ and $\tau_2w_0$ (one for the configuration, the
other for the velocity) are neutral for $S^{\Bbb R}x_0$, i. e. they do not
cause any change in the velocity history of the semi-trajectory.

\medskip

\proclaim{Lemma 4.4 (Lemma on the ``Neutral Trapezoid'')} Assume that 
$|\tau_i|<\epsilon_1$, $\tau_1\cdot\tau_2\ge0$, $t>0$. Introduce the notations
$t^*=(1+\tau_2^2)^{-1/2}t$, $x^*=S^{t^*}x_0=(q_0^*,v_0^*)$,
$(w_0^*,\,0)=\left(DS^{t^*}\right)(w_0,\,0)$, $\tau_1^*=\tau_1+t^*\tau_2$. 
Then the phase point $S^tx$ ($x=T_{\tau_1,\tau_2}x_0$) is equal to

$$
T_{\tau_1^*,\tau_2}x^*=\left(q_0^*+\tau_1^*w_0^*,\,(1+\tau_2^2)^{-1/2}
(v_0^*+\tau_2w_0^*)\right),
$$
provided that the so called ``neutral trapezoid''
$$
NT(x_0,w_0,\tau_1,\tau_2,t)=\left\{S^{t'}\left(T_{\tau'_1,\tau'_2}x_0
\right)\colon\; |\tau'_i|\le|\tau_i|,\; \tau'_i\cdot\tau_i\ge0,\;
0\le t'\le t\right\}
\tag 4.5
$$
is free of singularities, i. e. of multiple or tangential collisions.
\endproclaim

\medskip

\subheading{Remark 4.6} It should be noted, however, that --- when defining 
the translates
$$
T_{s,\tau_2}x^*=\left(q_0^*+s\cdot w_0^*,\, (1+\tau_2^2)^{-1/2}
(v_0^*+\tau_2w_0^*)\right)
$$
($0\le s\le\tau_1^*$) --- we may hit the boundary of the phase space, 
which means
that the time moment of a collision reaches the value zero. In that case,
in order to continue these translations beyond $s$,
by definition, we reflect both the direction vector $w_0^*$ of the spatial
variation and the velocity $(1+\tau_2^2)^{-1/2}(v_0^*+\tau_2w_0^*)$
with respect to the tangent hyperplane of the boundary $\partial\bold Q$
of the configuration space at the considered point of reflection. Although
we use this reflection in our definition of $T_{s,\tau_2}x^*$, in order to
keep our notations simpler, we do not indicate the arising change in
$w_0^*$ and $v_0^*+\tau_2w_0^*$ in the formulas. This convention will not
cause any confusion in the future. \endrem

\medskip

\subheading{Proof of Lemma 4.4} The whole point is that --- as long as we do
not hit any singularity --- there is a neutral $\Bbb R^2$-action
$A_{\alpha,\beta}$ ($(\alpha,\beta)\in\Bbb R^2$) lurking in the background:

$$
\cases
A_{\alpha,\beta}x_0=S^\alpha\left(T_{\beta,0}x_0\right), \\
A_{\alpha,\beta}\left(A_{\alpha',\beta'}x_0\right)=
S^{\alpha+\alpha'}\left(T_{\beta+\beta',0}x_0\right)
\endcases
\tag 4.7
$$
acting on the sheet 
$\left\{A_{\alpha,\beta}x_0\colon\;(\alpha,\beta)\in\Bbb R^2\right\}$,
where the phrase ``neutral'' means that there is no change in the velocity
process (i. e. no change in the collision normal vectors) during this 
perturbation. There is no problem with the neutrality 
and smoothness of this action $A_{\alpha,\beta}$ as long as we know that the
rectangle 

$$
\left\{A_{\alpha',\beta'}x_0\colon\;|\alpha'|\le|\alpha|,\; |\beta'|\le|\beta|,
\alpha\cdot\alpha'\ge0,\; \beta\cdot\beta'\ge0\right\}
$$
does not hit any singularity. The neutral trapezoid
$NT=NT(x_0,w_0,\tau_1,\tau_2,t)$ of (4.5) can be expressed in terms of the
$A_{\alpha,\beta}$-action as follows:

$$
\aligned
& NT(x_0,w_0,\tau_1,\tau_2,t)
=\big\{T_{0,\lambda}A_{\alpha,\beta}x_0\colon\; 0\le\alpha\le t^*, \\
& |\beta|\le|\tau_1|+\alpha|\tau_2|,\; |\lambda|\le|\tau_2|,\; 
\lambda\cdot\tau_i\ge0,\; \beta\cdot\tau_i\ge0\big\}.
\endaligned
\tag 4.8
$$
The statement of Lemma 4.4 is then easily provided by the commutativity and
neutrality of the $A_{\alpha,\beta}$-action. \qed

\medskip

Let $b>0$ be a suitably big number so that
$$
\text{dim}\Cal N_0\left(S^{[0,b]}x_0\right)=
\text{dim}\Cal N_0\left(S^{[0,\infty)}x_0\right).
\tag 4.9
$$
(The number $b$ is assumed to be not a moment of collision.) By further
strengthening (4.1), we may assume that the threshold $b>0$ is already
chosen so big and the radius $\epsilon_0$ of the ball 
$B_0=B(x_0,\, \epsilon_0)$ is selected so small that for every phase point
$x\in J\cap B_0$
$$
\text{dim}\Cal N_0\left(S^{[0,b]}x\right)=
\text{dim}\Cal N_0\left(S^{[0,\infty)}x\right)=
\text{dim}\Cal N_0(S^{[0,\infty)}x_0),
\tag 4.10
$$
and, on the other hand, $\text{dim}\Cal N_0\left(S^{[0,b]}x\right)=1$
for all $x\in B_0\setminus J$. Then, by selecting $\epsilon_1>0$
sufficiently small, we may assume that for
$|\tau_i|<\epsilon_1$, $\tau_1\cdot\tau_2\ge0$, the translated phase point
$x=T_{\tau_1,\tau_2}x_0$ is in $J\cap B_0$.

\medskip

\proclaim{Key Lemma 4.11}
For a typically selected phase point $x_0\in J$ (more precisely, apart from
a first category subset of $J$) the following holds true:

For every pair of real numbers $(\tau_1,\tau_2)$ with $|\tau_i|<\epsilon_1$,
$\tau_1\cdot\tau_2\ge0$, the positive trajectory $S^{[0,\infty)}x$ of
$x=T_{\tau_1,\tau_2}x_0$ ($\in J\cap B_0$) does not hit any singularity.
\endproclaim

\medskip

\subheading{Proof} We will argue by the absurd. Suppose that 
$S^t\left(T_{\tau_1,\tau_2}x_0\right)$ hits a singularity at time moment
$t=t_0$ ($>0$). Due to the smoothness of the orbit segments 
$S^{[0,b]}\left(T_{\tau_1,\tau_2}x_0\right)$, the number $t_0$ is necessarily
greater than $b$. The considered singularity can be one of the following
two types:

\medskip

\subheading{Type I. Tangential collision between the disks $i$ and $j$ at
time $t=t_0$} To simplify the notations, we assume that $\tau_i\ge0$,
$i=1,2$. Let us understand the relationship between the curve
$\gamma(s)=T_{\tau_1+s,\tau_2}x_0$ ($|s|<<1$) and the semi-invariant hull
$\bigcup_{t<0}S^t(\Cal S)$ of the considered tangential singularity
$\Cal S$. Due to the doubly neutral nature of the perturbations
$T_{\tau_1,\tau_2}x_0$, for the parameter values $s<0$ the disks
$i$ and $j$ must avoid each other (pass by each other) for 
$t\approx t_0$. Otherwise, if these disks collided on the orbit of 
$\gamma(s)$ near $t=t_0$ for $s<0$, then the further neutral perturbations
$\gamma(s)$ with $s\nearrow 0$ would not set these disks apart near $t=t_0$,
due to the neutral nature of the perturbations $T_{\tau_1+s,\tau_2}x_0$. 
Thanks to the neutrality of the perturbations 
$\gamma(s)=T_{\tau_1+s,\tau_2}x_0$ ($|s|<<1$), the smallest distance $d(s)$
($s<0$) between the disks $i$ and $j$ flying by each other around the time 
$t\approx t_0$ is an (inhomogeneous) linear function of the perturbation
parameter $s$ for $s<0$. Since, according to our assumption, for the value
$s=0$ the curve $\gamma(s)$ hits $\bigcup_{t<0}S^t(\Cal S)$, we get that
the constant derivative $d'(s)$ has to be negative for $s<0$, thus the curve
$\gamma(s)$ is transversal to the manifold $\bigcup_{t<0}S^t(\Cal S)$ at
$\gamma(0)$. Since further perturbations of $\gamma(s)$ with $s>0$ cause the
normal vector of the arising $(i,j)$-collision (at time $t\approx t_0$) to 
rotate, the perturbation direction vector $w_0$ turns out to be no longer
neutral with respect to the new collision. Thus, the above mentioned
transversality ``kills'' the neutral vector
$w_0\in\Cal N_0(S^{[0,\infty)}x)=\Cal N_0(S^{[0,b]}x)$
($x=T_{\tau_1,\tau_2}x_0=\gamma(0)$) by lowering the dimension of
$\Cal N_0\left(S^{[0,\infty)}\gamma(s)\right)$ for $s>0$. (We note that
new neutral vectors cannot appear because of the stable nature of 
$\Cal N_0\left(S^{[0,b]}\gamma(s)\right)$.) The latest statement, however,
contradicts to the assumed minimality of 
$\text{dim}\Cal N_0\left(S^{[0,\infty)}x\right)=
\text{dim}\Cal N_0\left(S^{[0,\infty)}x_0\right)$, see also (4.1) and the text
surrounding it.

\medskip

\subheading{Type II. A multiple collision singularity (of type 
$(i,j)$--$(j,k)$) at $t=t_0$}

\medskip

{\bf We begin with an important remark.} The multiple collision singularity 
of type $(i,j)$--$(j,k)$ means that on each side of the singularity a finite
sequence of alternating collisions $(i,j)$ and $(j,k)$ takes place in such a 
way that on one side of the singularity this finite sequence starts with
$(i,j)$, while on the other side it starts with $(j,k)$. This is how a
trajectory bifurcates into two different ``branches'', see \S\S2.3 above
about the notion of trajectory branches. Purely to simplify
the notations, hereby we are presenting a study of this type of singularity
in the case when both collision sequences are made up by two collisions.
This is not a restriction of generality, but merely a simplification of
the notations.

\medskip

Just as above, we again consider the curve
$\gamma(s)=T_{\tau_1+s,\tau_2}x_0$ ($|s|<<1$) and its relationship with the
invariant hull $\bigcup_{t\in\Bbb R}S^t(\Cal S)$ of the considered double
collision singularity $\Cal S$. Suppose that on the side $s<0$ of
$\bigcup_{t\in\Bbb R}S^t(\Cal S)$ the collision $(i,j)=\sigma_l$ precedes
the collision $(j,k)=\sigma_{l+1}$. Then the derivative of the time difference
$t(\sigma_{l+1})-t(\sigma_{l})$ with respect to $s$ (which depends on $s$
linearly, thanks to the neutrality of the vector $w_0$) at $s=0$ must be 
negative,
and, therefore, the curve $\gamma(s)$ transversally intersects the invariant 
hull $\bigcup_{t\in\Bbb R}S^t(\Cal S)$ of the studied double singularity at 
the point $\gamma(0)$. More precisely, denote by $v_i^-$, $v_j^-$, and
$v_k^-$ the velocities of the disks $i$, $j$, $k$ right before the collision
$\sigma_l$ on the trajectory of $\gamma(s)=T_{\tau_1+s,\tau_2}x_0$ for
$s<0$, $|s|<<1$. Similarly, let $v_i^+$, $v_j^+$, and $v_k^+=v_k^-$ the
corresponding velocities between $\sigma_l$ and $\sigma_{l+1}$ on the orbit
of the phase point $\gamma(s)$ for $s<0$, $|s|<<1$, and let 
$$
\cases
w_0^-=(\delta q_1^-,\dots,\delta q_N^-)=
\left(DS^{t(\sigma_{l})-\epsilon}(\gamma(s))\right)(w_0), \\
w_0^+=(\delta q_1^+,\dots,\delta q_N^+)=
\left(DS^{t(\sigma_{l})+\epsilon}(\gamma(s))\right)(w_0),
\endcases
\tag 4.12
$$
for $s<0$, $|s|<<1$. By the neutrality of $w_0^-$ and by the conservation
of the momentum we immediately obtain
$$
\cases
\delta q_i^--\delta q_j^-=\alpha(v_i^--v_j^-), \\
\delta q_i^+-\delta q_j^+=\alpha(v_i^+-v_j^+), \\
\delta q_i^+-\delta q_i^-=\alpha(v_i^+-v_i^-), \\
\delta q_k^+=\delta q_k^-,
\endcases
\tag 4.13
$$
where $\alpha$ is the advance of the collision $\sigma_l=(i,j)$ with respect
to the neutral vector $w_0$, see also \S 2. From the equations (4.13)
and from the conservation of the momentum (which is obviously also true for
the components $\delta q_a$ of neutral vectors) we obtain
$$
\cases
\delta q_i^+=\delta q_i^-+\alpha(v_i^+-v_i^-), \\
\delta q_j^+=\delta q_j^-+\alpha(v_j^+-v_j^-), \\
\delta q_j^+-\delta q_k^+=\beta(v_j^+-v_k^+)
=\delta q_j^--\delta q_k^-+\alpha(v_j^+-v_j^-).
\endcases
\tag 4.14
$$
Here $\beta$ denotes the advance of the collision $\sigma_{l+1}=(j,k)$ with
respect to $w_0$.

Let us study now the quite similar phenomenon on the other side of the
singularity $\bigcup_{t\in\Bbb R}S^t(\Cal S)$, i. e. for $s>0$. Since
$$
\Cal N_0\left(S^{[0,b]}\gamma(s)\right)=
\Cal N_0\left(S^{[0,b]}x_0\right)=
\Cal N_0\left(S^{[0,\infty)}x_0\right)
$$
and $\Cal N_0\left(S^{[0,\infty)}x_0\right)$ has the minimum value of all 
such dimensions, we obtain that 
$w_0\in\Cal N_0\left(S^{[0,\infty)}\gamma(s)\right)$. Thus, similar thing can 
be stated about the velocities and neutral vectors as above. Namely, denote by
$\tilde v_j^+$, $\tilde v_k^+$, and $\tilde v_i^+=v_i^-$ the velocities of
the disks $j$, $k$, $i$ between the collisions $\sigma_l=(j,k)$ and
$\sigma_{l+1}=(i,j)$ (Observe that the order of the two collisions is now
inverted!) on the orbit $S^{[0,\infty)}\gamma(s)$, $s>0$, $s<<1$. Let,
moreover, 
$$
\tilde w_0^+=(\delta\tilde q_1^+,\dots,\delta\tilde q_N^+)=
\left(DS^{t(\sigma_{l})+\epsilon}(\gamma(s))\right)(w_0)
\tag 4.15
$$
for $s>0$, $s<<1$, and $\tilde\beta$, $\tilde\alpha$ be the advances of the
collisions $\sigma_l=(j,k)$, and $\sigma_{l+1}=(i,j)$, respectively. Then
$\delta q_j^--\delta q_k^-=\tilde\beta(v_j^--v_k^-)$ in the last equation of
(4.14), so we get that
$$
\beta(v_j^+-v_k^+)=\tilde\beta(v_j^--v_k^-)+\alpha(v_j^+-v_j^-).
\tag 4.16
$$
By neutrality, for all orbits $S^{[0,\infty)}\gamma(s)$ ($|s|<<1$) the
$(i,j)$ collision near $t=t_0$ (which is either $\sigma_l$ or 
$\sigma_{l+1}$, depending on which side of the singularity we are) has the 
same normal vector $\vec n_1$ and, similarly, for all orbits 
$S^{[0,\infty)}\gamma(s)$ the $(j,k)$ collision near $t=t_0$ has the same
normal vector $\vec n_2$. How can we take now advantage of (4.16)? 
First of all, we can assume that the relative velocities 
$v_j^--v_k^-$ and $v_j^+-v_k^+$ are nonzero, for each of the equations
$v_j^--v_k^-=0$ and $v_j^+-v_k^+=0$ defines
a codimension-two set, which is atypical in $J$, so we can assume that these
vectors are nonzero on the orbit of $\gamma(0)$ (or, equivalently, on the
orbit of $x_0$) by typically choosing the starting phase point $x_0$.
Secondly, by adding an appropriate scalar multiple of $v_0$ to the neutral
vector $w_0$, we can achieve that $\alpha=0$ in (4.16), see also \S 2.
We infer, therefore, that the relative velocities
$v_j^--v_k^-$ and $v_j^+-v_k^+$ are parallel, as long as at least one of the
advances $\beta$ and $\tilde\beta$ is nonzero. However, 
$\alpha=\beta=\tilde\beta=0$ would mean that 
$\delta q_i^-=\delta q_j^-=\delta q_k^-$, which is impossible, for in that
case the time difference $t(\sigma_{l+1})-t(\sigma_{l})$ would not change
(and, therefore, it could not tend to zero) as $s\nearrow0$. Thus, we conclude
that $v_j^--v_k^-\parallel v_j^+-v_k^+$. However, the difference of these 
vectors is obviously parallel to the collision normal $\vec n_1$, so we get
$$
v_j^--v_k^-\parallel\vec n_1.
\tag 4.17
$$
A similar argument yields
$$
v_i^--v_j^-\parallel\vec n_2.
\tag 4.18
$$
However, the events described in (4.17--18) together define a codimension-two
subset of the phase space, so we can assume that the typically selected
starting phase point $x_0\in J$ is outside of all such codimension-two 
sub-manifolds. This finishes the proof of Main Lemma 4.11. \qed

\bigskip \bigskip

\heading
\S 5. Non-Existence of $J$-Manifolds. \\
Part B: The Weird Behavior of the $\Omega$-limit Set
\endheading

\bigskip \bigskip

Let us study now the non-empty, compact $\Omega$-limit set
$$
\Omega(x_0)=\Big\{x_\infty\in\bold M\colon\; \exists\text{ a sequence }
t_n\nearrow\infty\text{ such that } x_\infty=\lim_{n\to\infty}S^{t_n}x_0\Big\}
\tag 5.1
$$
of the trajectory $S^{\Bbb R}x_0=\left\{S^tx_0=x_t\colon\; t\in\Bbb R\right\}$.
Consider an arbitrary phase point $x_\infty\in\Omega(x_0)$,
$x_\infty=\lim_{n\to\infty}x_{t_n}$, $t_n\nearrow\infty$. Although the 
trajectory of $x_\infty$ may be singular, we can assume that we have properly
selected and fixed a branch $S^{(-\infty,\infty)}x_\infty$ of the trajectory
of $x_\infty$ (for the notion of trajectory branches, please see \S\S2.3
above), so that
whenever $x_\infty$ belongs to a singularity $S^{-t}\Cal S$, the
sequence of points $x_{t_n}$ converges to $x_\infty$ from one side of the
codimension-one sub-manifold $S^{-t}\Cal S$. This can be achieved by using
Cantor's diagonal method and switching to a subsequence of the sequence
$t_n\nearrow\infty$. Then for $t\in\Bbb R$ the phase points $x_{t_n+t}$
will converge as $n\to\infty$, and we will define the limit
$\lim_{n\to\infty}x_{t_n+t}$ as $S^tx_\infty$. In this way we correctly define
a trajectory branch $S^{\Bbb R}x_\infty$ of the phase point $x_\infty$. As for
the concept of trajectory branches, see \S2.3.

By switching again -- if necessary -- to a suitable subsequence of 
$t_n\nearrow\infty$, we can assume that the unit neutral vectors
$w_{t_n}=\left(DS^{t_n}(x_0)\right)(w_0)$ converge to a (unit) neutral
vector $w_\infty\in\Cal N_0(S^{\Bbb R}x_\infty)$, which is then necessarily
perpendicular to the velocity 
$v_\infty=v(x_\infty)=(v_1^\infty,\dots,v_N^\infty)$ of $x_\infty$. We write
$$
w_\infty=\left(\delta q_1^\infty,\dots,\delta q_N^\infty\right),\quad
x_\infty=\left(q_1^\infty,\dots,q_N^\infty;\,v_1^\infty,
\dots,v_N^\infty\right).
$$
We would like to point out again that the well defined orbit 
$S^{\Bbb R}x_\infty$ may be singular. In the case of a multiple
collision, according to what was said above, the infinitesimal time-ordering
of the collisions (taking place at the same time) is determined, just as the 
resulting product of reflections connecting the incoming velocity $v^-$ with
the outgoing velocity $v^+$. As far as the other type of singularity --- the
tangential collisions --- is concerned, here there are two possibilities.
The first one, in which case the tangentially colliding disks $i$ and $j$ 
have proper collisions on the nearby approximating trajectories
$S^{\Bbb R}x_{t_n}$, $n\to\infty$. The second case is when the
tangentially colliding disks $i$ and $j$ pass by each other without collision
on the approximating orbit $S^{\Bbb R}x_{t_n}$, $n\to\infty$. In both cases, 
we do not include a tangential collision in the symbolic collision 
sequence of $S^{\Bbb R}x_{\infty}$. In the sequel we will exclusively
deal with non-tangential collisions, i. e. collisions with nonzero momentum
exchange. They are called proper collisions. This note has particular 
implications when defining the connected components of the collision graph
of the entire trajectory $S^{\Bbb R}x_{\infty}$.

\subheading{Definition 5.2} Let
$\{1,2,\dots,N\}=H_1\cup H_2\cup\dots\cup H_k$ be the partition of the vertex
set into the connected components of the collision graph 
$\Cal G\left(S^{\Bbb R}x_{\infty}\right)$ of the orbit $S^{\Bbb R}x_{\infty}$.
For any $i$, $1\le i\le k$, we denote by
$\{S_i^t\}_{t\in\Bbb R}=\{S_i^t\}$ the internal dynamics of the subsystem
$H_i$, i. e. the dynamics in which we

(a) reduce the total momentum of the subsystem $H_i$ to zero by observing it
from a suitably moving reference system;

(b) do not make any distinction between two configurations of $H_i$ differing
only by a uniform spatial translation; (Factorizing with respect to uniform 
spatial translations, see also \S 1.)

(c) carry out a time-rescaling, so that the total kinetic energy of the 
internal system $\{S_i^t\}_{t\in\Bbb R}$ is equal to $1$. \endrem

\medskip

Let, moreover, $M_i=\sum_{j\in H_i}m_j$ the total mass, 
$I_i=\sum_{j\in H_i}m_jv_j^\infty$ the total momentum, and $V_i=I_i/M_i$
the average velocity of the subsystem $H_i$. Similarly, we write 
$W_i=(M_i)^{-1}\sum_{j\in H_i}m_j\delta q_j^\infty=
(M_i)^{-1}\sum_{j\in H_i}m_jw_j^\infty$ for the total (average) displacement
of the system $H_i$ under the action of the neutral vector 
$w_\infty=(\delta q_1^\infty,\dots,\delta q_N^\infty)=
(w_1^\infty,\dots,w_N^\infty)$. Finally, let $|H_i|\ge2$ for 
$i\in\{1,2,\dots,s\}$, $|H_i|=1$ for $i\in\{s+1,\dots,k\}$. 

\medskip

First of all, we prove

\proclaim{Lemma 5.3} Let $\lambda,\, \mu\in\Bbb R$ be given numbers, and
form the neutral vector 
$$
n(\mu,\lambda)=\mu v_\infty+\lambda w_\infty\in\Cal N_0(x_\infty)=
\Cal N_0(S^{\Bbb R}x_\infty).
$$
Define 
$T_{n(\mu,\lambda),0}x_\infty=\left(q_\infty+n(\mu,\lambda),v_\infty\right)$
as the neutral translation of
$x_\infty=\mathbreak(q_\infty,v_\infty)$ by the vector 
$n(\mu,\lambda)$, where we use the natural convention of Remark 4.6. Let,
finally, $i$ and $j$ be labels of disks belonging to different components
$H_l$, say, to $H_p$ and $H_q$. We claim that the orbit 
$S^{\Bbb R}T_{n(\mu,\lambda),0}x_\infty$ of $T_{n(\mu,\lambda),0}x_\infty$
cannot have a proper (i. e. non-tangential) collision between the disks
$i$ and $j$.
\endproclaim

\subheading{Proof} Assume the contrary. Then, by a simple continuity argument,
one finds some real numbers $\mu_0$, $\lambda_0$ for which the orbit of
$T_{n(\mu_0,\lambda_0),0}x_\infty$ hits a tangential singularity between the
disks $i$ and $j$. By using a suitably accurate approximation 
$(x_{t_n},w_{t_n})\approx(x_\infty,w_\infty)$, one finds a neutral, spatial
translation of $x_{t_n}$ by a vector $\mu_1v_{t_n}+\lambda_1w_{t_n}$
($\mu_1\approx\mu_0$, $\lambda_1\approx\lambda_0$) such that the orbit of
$(q_{t_n}+\mu_1v_{t_n}+\lambda_1w_{t_n},\, v_{t_n})$ hits a tangential
singularity between the disks $i$ and $j$, which is impossible by Lemma 4.11.
This finishes the indirect proof of 5.3. \qed

\medskip

The main step in the indirect proof of the Theorem is

\proclaim{Key Lemma 5.4}
There exists a finite collection of nonzero lattice vectors 
$$
l_0,l_1,\dots,l_p\in\Bbb Z^2
$$
(depending only on $N$ and the common
radius $r$ of the $N$ disks moving in the standard unit torus 
$\Bbb T^2=\Bbb R^2/\Bbb Z^2$) with the following properties: For every 
separating manifold $J$, for every phase point $x_0\in J$ fulfilling 
conditions (0)---(3) listed at the end of \S3, and for
every $\Omega$-limit point $x_\infty=\lim_{n\to\infty}x_{t_n}$ of the orbit
$S^{\Bbb R}x_0$ it is true that $k\ge2$ (i. e. the collision graph $\Cal G$
of $S^{\Bbb R}x_\infty$ is not connected, see Definition 5.2 above),
and there is an index $j\in\{0,1,\dots,p\}$ such that all velocities
$v_i^\infty$ ($i=1,\dots,N$) are parallel to the lattice vector $l_j$.
\endproclaim

\subheading{Remark 5.4/a} It is easy to see that the scenario described in the
key lemma (i. e. that all velocities are parallel to $l_j$ for all time 
$t\in\Bbb R$) can only take place if the dynamically connected components of
the motion --- the connected components of the collision graph of
$S^{\Bbb R}x_\infty$ --- move on closed geodesics of $\Bbb T^2$ being parallel
to $l_j$. \endrem

\medskip

\subheading{Remark 5.4/b} The part $k\ge2$ of the key lemma does not play
any role in the overall proof of the Main Theorem. The reason why we included
it is of didactics: When indirectly proving $k\ge2$ below (under section 
$1^o$) we obtain an auxiliary result saying that the advances of a connected 
subsystem are necessarily equal, and this will be later used in proving the 
key lemma for the general case $k\ge2$.

\medskip

\subheading{Proof of Key Lemma 5.4} First of all, we prove the geometric

\proclaim{Sub-lemma 5.5}
Consider the standard $x$---$y$ coordinate plane with the
usual unit vectors $e_1=(1,0)$, $e_2=(0,1)$. Suppose that infinitely many
disks of radius $r$ and centers at $q_i+je_2\in\Bbb R^2$ ($i=1,\dots,N$;
$j\in\Bbb Z$) are moving uniformly in $\Bbb R^2$ and colliding elastically.
We assume that the disk centered at $q_i+je_2$ has mass $m_i$, and its velocity
$v_i=\dot q_i$ is also independent of $j$, $i=1,\dots,N$, $j\in\Bbb Z$. We 
claim that if the trajectory of such an $e_2$-periodic system remains in the
half plane $x\le L$ (for all time $t\in\Bbb R$, the number $L$
is given), then all velocities $\dot q_i(t)=v_i(t)$ are parallel to $e_2$.
\endproclaim

\subheading{Proof of 5.5} We carry out an induction on the number of disks
$N$ of the $e_2$-factorized system. For $N=1$ the statement is obviously true.
Let $N>1$, and assume that the sub-lemma has been proved for all numbers
$N'<N$. 

Let $i_1,i_2,\dots,i_a$ be the labels $i$ of disks with the largest value
of the inner product $\langle q_i,e_1\rangle$ at time $t=0$. To simplify the
notations, we assume that the indices $i_1,i_2,\dots,i_a$ are just
$1,2,\dots,a$. 

Suppose first that the $x$-coordinate $\langle v_i,e_1\rangle$ of the velocity
$v_i=\dot q_i$ is nonzero for some $i\le a$. By reversing time, if necessary,
we can assume that $\langle v_i,e_1\rangle>0$ for some $i\le a$. This means
that among the disks with the rightmost position at least one moves to the
right. Denote by $i=1,2,\dots,b$ ($1\le b\le a$) the labels of disks $i$
($i\le a$) for which the inner product $\langle v_i,e_1\rangle$ is maximal.
Now it is easy to see that the first velocity component
$\langle v_1,e_1\rangle=\dots=\langle v_b,e_1\rangle$ ($>0$) cannot decrease
in time. As a matter of fact, two things can only happen to the disk(s) $i$ 
with the rightmost position and maximum value of $\langle v_i,e_1\rangle$:
Either the disk $i$ collides with another disk coming from the left, or
another disk $j$ with a larger velocity component $\langle v_j,e_1\rangle$
passes by $q_i$, thus by snapping the ``title'' of having the rightmost 
position. In either case, the maximum value of the first velocity component
$\langle v_i,e_1\rangle$ of the rightmost disk(s) can only increase. This
argument shows that at least one disk $i$ will escape to the right
($\langle q_i,e_1\rangle\to+\infty$), which is impossible by our assumption
on the boundedness of the $x$-coordinates. 

Therefore, only the second possibility can occur, i. e. that 
$\langle v_i,e_1\rangle=0$ for all $i\le a$. This should then remain valid
for all time $t\in\Bbb R$ by the above argument. However, this also means that
the disks $q_i+je_2$, $i\le a$, $j\in\Bbb Z$, collide among themselves, while
all of them have vertical velocities. In the case $a=N$ we are done, while in
the case $a<N$ we can use the induction hypothesis, which says that all
velocities $v_i$, $i>a$, are also vertical. This finishes the proof of
5.5. \qed

\medskip

\proclaim{Remarks 5.5/a} 

1. By taking a brief look at the proof, we can see that it readily generalizes
to any dimension $d\ge2$. What is even more, if the single boundedness
condition of the lemma is replaced by $k$ linearly independent linear
inequalities $A_j\left(q_i(t)\right)\le L_j$ (for all $t\in\Bbb R$, 
$i=1,\dots,N$, $j=1,\dots,k$, the linear functionals $A_j$ being linearly 
independent), then we can state that all velocities $v_i(t)$ ($t\in\Bbb R$,
$i=1,\dots,N$) belong to some $d-k$-dimensional subspace $S$ of $\Bbb R^d$,
and the positions line up in groups on translated copies of the subspace $S$.

2. The postulated periodicity ($e_2$--periodicity) condition has not been used
in the proof and, therefore, it can be dropped.
\endproclaim

\medskip

Let us return now to the proof of Key Lemma 5.4. Its proof will be divided
into several parts.

\medskip

\subheading{$1^o$ First we prove that $k\ge2$, i. e. the collision graph
$\Cal G$ of $S^{\Bbb R}x_\infty$ is not connected} Assume, on the contrary,
the connectedness of $\Cal G$. Let us focus on the limiting neutral vector
$$
w_\infty=(\delta q_1^\infty,\dots,\delta q_N^\infty)
=(w_1^\infty,\dots,w_N^\infty)=\lim_{n\to\infty}w_{t_n},
$$ 
for which $\sum_{i=1}^N m_iw_i^\infty=0$, 
$\sum_{i=1}^N m_i||w_i^\infty||^2=1$, and
$\sum_{i=1}^N m_i\langle w_i^\infty,v_i^\infty\rangle=0$, where
$x_\infty=(q_1^\infty,\dots,q_N^\infty;\, v_1^\infty,\dots,v_N^\infty)$.
Let $\Sigma=(\dots,\sigma_{-1},\sigma_0,\sigma_1,\dots)$ be the symbolic
collision sequence of $S^{\Bbb R}x_\infty$, and denote by 
$\alpha_j=\alpha(\sigma_j)$ the advance of the collision $\sigma_j$ with 
respect to the neutral vector $w_\infty$, see \S 2. Since $w_\infty$
is not parallel to $v_\infty=(v_1^\infty,\dots,v_N^\infty)$, by using the 
assumed connectedness of $\Cal G$ we get that not all advances $\alpha_j$
($j\in\Bbb Z$) are equal, see the second statement of Lemma 2.13 in
[Sim(1992-B)]. (That statement says that, in the case of a connected collision
graph $\Cal G$, the equality of all advances $\alpha_j$ implies that the 
considered neutral vector $w_\infty$ is parallel to the velocity $v_\infty$.)
By switching from $w_\infty$ to $-w_\infty$, if necessary, we
can assume that there are indices $j<k$ ($j,\, k\in\Bbb Z$) for which 
$\alpha_j<\alpha_k$, $\sigma_j\ne\sigma_k$, and 
$\sigma_j\cap\sigma_k\ne\emptyset$.
This means, however, that the translated copy
$(q_\infty+\lambda w_\infty, v_\infty)$ of $x_\infty=(q_\infty, v_\infty)$
will hit a double collision singularity $t(\sigma_l)=t(\sigma_{l+1})$
($\sigma_l\cap\sigma_{l+1}\ne\emptyset$) for some value
$$
0<\lambda\le\lambda^*=(\alpha_{k}-\alpha_{j})^{-1}\cdot
\left(t(\sigma_{k})-t(\sigma_{j})\right).
$$
By considering some well approximating pair 
$(x_{t_n},w_{t_n})\approx(x_\infty,w_\infty)$, we get that some translated
copy $(q_{t_n}+\lambda'w_{t_n}, v_{t_n})$ of $x_{t_n}=(q_{t_n},v_{t_n})$
will hit a double collision singularity for some $\lambda'\approx\lambda$.
However, this statement contradicts to Lemma 4.11. Therefore, the collision 
graph $\Cal G$ of $S^{\Bbb R}x_\infty$ is not connected, the number $k$ of
the connected components of $\Cal G$ is at least two. \qed

\medskip

\subheading{$2^o$ Next we prove Key Lemma 5.4 in the case $k=N$, i. e. when
no proper collision at all takes place on the trajectory $S^{\Bbb R}x_\infty$}
(By the way, this phenomenon can only occur if the maximum lengths
$\tau(x_{t_n})$ and $\tau(-x_{t_n})$ of the collision-free paths of
$x_{t_n}=(q_{t_n},v_{t_n})$ and $-x_{t_n}=(q_{t_n},-v_{t_n})$ tend to
infinity, as $n\to\infty$.)

\medskip

First we put forward

\proclaim{Sub-lemma 5.6} For every pair of labels $(i,j)$ ($1\le i<j\le N$)
it is true that
$$
\dim\text{span}\{v_i^\infty-v_j^\infty,\, w_i^\infty-w_j^\infty\}\le 1.
$$
\endproclaim

\medskip

This sub-lemma is an immediate consequence of Lemma 5.3.

\medskip

By the normalizations $\sum_{i=1}^N m_iv_i^\infty=0$ and 
$\sum_{i=1}^N m_i||v_i^\infty||^2=1$, not all velocities 
$v_1^\infty,\dots,v_N^\infty$ are the same.

\proclaim{Sub-lemma 5.7} All velocities $v_1^\infty,\dots,v_N^\infty$ are
parallel to each other, that is, 

$$
\dim\text{span}\{v_1^\infty,\dots,v_N^\infty\}=1.
$$

\endproclaim

\subheading{Proof} Assume the opposite, i. e. 
$\dim\text{span}\{v_1^\infty,\dots,v_N^\infty\}=2$. (Again an essential use of
the condition $\nu=2$.) Due to the relation
$\sum_{i=1}^N m_iv_i^\infty=0$, the points $v_1^\infty,\dots,v_N^\infty$
of the plane $\Bbb R^2$ do not lie on the same line, not even on a line not 
passing through the origin. Thus, the points
$v_1^\infty,\dots,v_N^\infty$ of $\Bbb R^2$ do not lie on the same affine
line. We may assume that $v_1^\infty$, $v_2^\infty$, and $v_3^\infty$
do not lie on the same affine line of $\Bbb R^2$. Since 
$w_2^\infty-w_1^\infty=\alpha(v_2^\infty-v_1^\infty)$,
$w_3^\infty-w_2^\infty=\beta(v_3^\infty-v_2^\infty)$, and
$w_3^\infty-w_1^\infty=\gamma(v_3^\infty-v_1^\infty)$, we conclude that

$$
\gamma(v_2^\infty-v_1^\infty)+\gamma(v_3^\infty-v_2^\infty)=
\alpha(v_2^\infty-v_1^\infty)+\beta(v_3^\infty-v_2^\infty),
$$
so $\alpha=\beta=\gamma$ by the linear independence of 
$v_2^\infty-v_1^\infty$ and $v_3^\infty-v_2^\infty$. For any index $i>3$
with $v_i^\infty\ne v_1^\infty$ we have that 
$(v_1^\infty,\,v_2^\infty,\,v_i^\infty,)$ or 
$(v_1^\infty,\,v_3^\infty,\,v_i^\infty,)$ do not lie on the same affine line, 
and again conclude (the same way as above) that 

$$
w_i^\infty-w_1^\infty=\alpha(v_i^\infty-v_1^\infty)
\tag 5.8
$$
with the same $\alpha$ as above. It is obvious that (5.8) also holds for
$i>3$ with $v_i^\infty=v_1^\infty$ and for $i=1,\,2,\,3$, i. e. (5.8) is true
for all $i=1,\dots,N$. Thanks to the conventions 
$\sum_{i=1}^N m_iw_i^\infty=\sum_{i=1}^N m_iv_i^\infty=0$, the equations 
(5.8) can only be fulfilled by the solution $w_i^\infty=\alpha v_i^\infty$
($i=1,\dots,N$), which is impossible, for the vector $w_\infty$ is not
parallel to $v_\infty$. This contradiction finishes the indirect proof of
Sub-lemma 5.7. \qed

\medskip

Now continue the proof of Key Lemma 5.4 in the case $k=N$. We got that all
velocities $v_i^\infty$ in $S^{\Bbb R}x_\infty$ are parallel to the same
direction vector $0\ne l\in\Bbb R^2$. Since the uniformly moving disks of
the orbit $S^{\Bbb R}x_\infty$ have no proper collision, we get that
$$
\text{dist}\{q_2^\infty-q_1^\infty+t\cdot l,\, 0\}\ge 2r
$$
for all $t\in\Bbb R$. This means, however, that the direction vector $l$ is
parallel to an irreducible (non-divisible) lattice vector 
$0\ne l_0\in\Bbb Z^2$, such that $||l_0||\le\dfrac{1}{4r}$. There are only 
finitely many choices for such a lattice vector $l_0\in\Bbb Z^2$. This
completes the proof of Key Lemma 5.4 in the case $k=N$. \qed

\medskip

\subheading{$3^o$ The case $s=k$ ($\ge 2$), i. e. when $|H_i|\ge2$ for all
$i$, $i=1,\dots,k$}

Let us study, first of all, the relationship between the subsystems $H_1$
and $H_2$ (and their internal dynamics $\{S_1^t\}$, $\{S_2^t\}$) with 
particular emphasis on their relation to the limiting neutral vector
$w_\infty\in\Cal N_0(S^{\Bbb R}x_\infty)$. Lemma 2.13 of [Sim(1992-B)] yields
(see also the reference to that result in the exposition of $1^o$ above)
that the advances of all collisions of $\{S_1^t\}_{t\in\Bbb R}$ with
respect to $w_\infty$ are equal to the same number $\alpha$ and, similarly,
all collisions of the internal flow $\{S_2^t\}_{t\in\Bbb R}$ share the same
advance $\beta$ with respect to the neutral vector $w_\infty$. Select and
fix an arbitrary real number $t_0$, and consider the linear combination
$$
n(\lambda)=(t_0-\alpha\lambda)v_\infty+\lambda w_\infty\in
\Cal N_0(S^{\Bbb R}x_\infty)
$$
with variable $\lambda\in\Bbb R$. Also consider the corresponding neutral
spatial translation 
$$
x_\infty=(q_\infty,\, v_\infty)\longmapsto T_{n(\lambda),0}x_\infty=
(q_\infty+n(\lambda),\, v_\infty)
$$
of $x_\infty$ with the natural convention of Remark 4.6. Observe that the 
neutral translation $T_{n(\lambda),0}$ has the following effect on the
internal dynamics $\{S_1^t\}_{t\in\Bbb R}$ and $\{S_2^t\}_{t\in\Bbb R}$:
The advance of the subsystem $H_1$ is $t_0$, i. e. the internal time of
evolution of $H_1$ will be the fixed number $t_0$. On the other hand, the
advance of $H_2$ is obviously $t_0+\lambda(\beta-\alpha)$. We distinguish
between two, quite differently behaving situations:

\medskip

\subheading{Case (A): $\alpha\ne\beta$} The internal time of the subsystem
$H_2$ (under the translation $T_{n(\lambda),0}x_\infty$, now 
$\lambda\in\Bbb R$ plays the role of time) changes linearly with $\lambda$,
it is equal to $t_0+\lambda(\beta-\alpha)$, while the internal time of
$H_1$ is constantly $t_0$. How about the relative motion of the
non-interacting groups $H_1$ and $H_2$? Recall that 
$V_i=(M_i)^{-1}\sum_{j\in H_i}m_jv_j^\infty$ is the average velocity of the
subsystem $H_i$, while $W_i=(M_i)^{-1}\sum_{j\in H_i}m_jw_j^\infty$ is the 
average displacement of the subsystem $H_i$ under the translation by the
neutral vector $w_\infty\in\Cal N_0(S^{\Bbb R}x_\infty)$. (Note that
$M_i=\sum_{j\in H_i}m_j$.) The relative position of the subsystem $H_1$ with
respect to $H_2$ can be measured, for example, by the relative position 
$q_{j_1}^\infty-q_{j_2}^\infty$ of the disks $j_1\in H_1$, $j_2\in H_2$,
$j_1$, $j_2$ fixed. To simplify the notations, we assume that $j_1=1$,
$j_2=2$. Thus the relative position of the subsystem $H_1$ with respect to
$H_2$ varies with $\lambda$ as follows:
$$
\aligned
q_1^\infty(\lambda)-q_2^\infty(\lambda)=q_1^\infty-q_2^\infty+
(t_0-\alpha\lambda)(V_1-V_2)+\lambda(W_1-W_2) \\
=q_1^\infty-q_2^\infty+t_0(V_1-V_2)+\lambda\left[W_1-W_2-\alpha(V_1-V_2)
\right].
\endaligned
\tag 5.9
$$
Now we would like to paint a global picture (global, that is, in the universal
covering space $\Bbb R^2$) of the orbit of $H_2$ under the neutral spatial
translations $T_{n(\lambda),0}$, $\lambda\in\Bbb R$. Due to the factorization
with respect to uniform spatial translations when defining our model 
(see \S 1), in order to lift the dynamics from $\Bbb T^2$ to its
universal covering space $\Bbb R^2$ (in a $\Bbb Z^2$-periodic manner),
it is necessary and sufficient to specify
the position of the lifted copy $\bar q_1^\infty(\lambda)\in\Bbb R^2$ of
$q_1^\infty(\lambda)=q_1\left(T_{n(\lambda),0}x_\infty\right)$. We take 
$$
\bar q_1(\lambda)=\bar q_1=\int_0^{t_0}\left(v_1(S^tx_\infty)-V_1\right)dt
$$
(independently of $\lambda$, so that the ``baricenter'' of $H_1$ is unchanged
while $t_0$ is changing later on), since the internal time of the subsystem
$H_1$ is constantly $t_0$, and we want to describe the motion of $H_2$
relative to $H_1$. For $i=1,2,\dots,N$ let the resulting $\Bbb Z^2$-periodic
lifting to 
$\Bbb R^2$ of $q_i^\infty(\lambda)=q_i\left(T_{n(\lambda),0}x_\infty\right)$ be
$$
\bar q_i(\lambda)+l, \qquad l\in\Bbb Z^2,
\tag 5.10
$$
where the lifting $\bar q_i(\lambda)\in\Bbb R^2$ is selected in such a way 
that it depends on $\lambda$ continuously. We point out here that currently
the translation parameter $\lambda$ plays the role of time. Also note that
for any $j\in H_1$ we have $\bar q_j(\lambda)=\text{const}$ (independent
of $\lambda$), for 
$\bar q_1(\lambda)=\int_0^{t_0}\left(v_1(S^tx_\infty)-V_1\right)dt$,
and the internal time of the subsystem $H_1$ is not changing by the 
translations $T_{n(\lambda),0}$. We want to pay special attention to the
orbit of points $\bar q_j(\lambda)+l\in\Bbb R^2$, $j\in H_2$, $l\in\Bbb Z^2$.
We define the open $2r$-neighborhood 
$U=U(x_\infty,w_\infty,H_1,H_2,t_0)$ of the set 
$$
\left\{\bar q_j(\lambda)+l\colon\;j\in H_2,\; l\in\Bbb Z^2,\; \lambda\in\Bbb R
\right\}
$$
as follows:
$$
U=\Big\{x\in\Bbb R^2\colon\; \exists\; j\in H_2,\; l\in\Bbb Z^2,\; 
\lambda\in\Bbb R
\text{ s. t. } \text{dist}\left(x,\, \bar q_j(\lambda)+l\right)<2r\Big\}.
\tag 5.11
$$
According to Lemma 5.3, the points $\bar q_j(\lambda)=\bar q_j$ ($j\in H_1$)
do not belong to the $\Bbb Z^2$-periodic open set $U$. Let us understand the
connected components of the set $U$. Since the open set $U$ is 
$\Bbb Z^2$-periodic, the $\Bbb Z^2$-translations will just permute the 
connected components of $U$ among themselves. The following lemma essentially
uses the $2-D$ topology of $\Bbb R^2$:

\medskip

\proclaim{Sub-lemma 5.12} Let $U_0\subset U$ be a connected component of a
$\Bbb Z^2$-periodic open set $U$. Then exactly one of the following
possibilities occurs:

(1) $U_0$ is bounded;

(2) $U_0$ is unbounded, $l_0$-periodic with some lattice vector 
$0\ne l_0\in\Bbb Z^2$, and $U_0$ is bounded in the direction perpendicular to
$l_0$;

(3) $U_0$ is $\Bbb Z^2$-periodic.
\endproclaim

\medskip

\subheading{Remark 5.13} In the case (2) all periodicity vectors $l\in\Bbb Z^2$
of $U_0$ are integer multiples of an irreducible lattice vector $l_0$, which
is uniquely determined up to a sign. \endrem

\subheading{Proof} Denote by $p:\;\Bbb R^2\to\Bbb T^2$ the natural projection.
Consider the open and connected set $V_0=p(U_0)\subset\Bbb T^2$. It follows
immediately from the conditions of the sub-lemma that

\medskip

(a) $U_0$ is a connected component of the open set $p^{-1}(V_0)$, and

\medskip

(b) $p:\;U_0\to V_0$ is a covering map.

\medskip

It is well known from the elements of topology that the group

$$
G=\left\{g\in\Bbb Z^2\big|\; U_0+g=U_0\right\}
$$
is the group of all deck automorphisms of the covering $p:\;U_0\to V_0$, and
$G$ is naturally isomorphic to the fundamental group of $V_0$. Now there are
three possibilities for the subgroup $G$ of $\Bbb Z^2$:

\medskip

(1) $|G|=1$;

(2) $G\cong\Bbb Z$ (i. e. $\text{rank}(G)=1$);

(3) $G\cong\Bbb Z^2$ (i. e. $\text{rank}(G)=2$).

\medskip

In the first case the covering $p:\;U_0\to V_0$ is an isometry, so $U_0$ is
bounded.

In the second case, let $l_0\in G$ be a generating element of $G$. Topology
says that the image $V_0$ of $U_0$ under the covering map $p:\;U_0\to V_0$
is just the factor of $U_0$ with respect to all translations by the integer
multiples of $l_0$. This means that (2) of 5.12 holds true.

Finally, in the third case the domain $U_0$ contains an $l_1$-periodic,
continuous curve $\gamma_1$ and an $l_2$-periodic, continuous curve $\gamma_2$,
where $l_1$ and $l_2$ are two linearly independent elements of $G$. 
It follows immediately from the
topology of the plane $\Bbb R^2$ that the translate $\gamma_1+(m,n)$
intersects the curve $\gamma_2$ for any $(m,n)\in\Bbb Z^2$. Since $U_0+(m,n)$
is a connected component of the $\Bbb Z^2$-periodic open set $U$ and 
$\left[U_0+(m,n)\right]\cap U_0\ne\emptyset$, we have that $U_0+(m,n)=U_0$
for any $(m,n)\in\Bbb Z^2$, and this is just case (3) of the sub-lemma. \qed

\medskip

The next sub-lemma takes into account that the open set $U_0$ (a connected
component of the open set $U$ defined in (5.11)) is determined
by a special dynamical system.

\proclaim{Sub-lemma 5.14} Out of the three cases listed in 5.12, in fact only
one of them, namely (2) can occur.
\endproclaim

\subheading{Proof} 

1. The impossibility of (1): Observe that for every $\bar q_j(\lambda)$
($j\in H_2$, $\lambda\in\Bbb R$) there exists a lattice vector $l\in\Bbb Z^2$
such that $\bar q_j(\lambda)+l\in U_0$. This follows simply from the
connectedness of the collision graph of the $H_2$ subsystem $\{S_2^t\}$.
If $U_0$ were bounded, then there would be a bounded cluster (enclosure)
of a billiard dynamics with positive kinetic energy inside $U_0$, which is
impossible for many reasons, for example, by Sub-lemma 5.5. Thus, $U_0$
is necessarily unbounded.

2. The impossibility of (3): Assume that $U_0$ is $\Bbb Z^2$-periodic.
Then $U_0$ contains an $e_1$-periodic, continuous curve $\gamma_1$, and
an $e_2$-periodic, continuous curve $\gamma_2$, as well. The 
$\Bbb Z^2$-periodic system of curves
$$
\bigcup\Sb m\in\Bbb Z\endSb(\gamma_1+me_2)\cup
\bigcup\Sb m\in\Bbb Z\endSb(\gamma_2+me_1)\subset U_0
$$
shows that the connected components of $\Bbb R^2\setminus U_0$ are bounded.
(Here we essentially use the $2-D$ topology of $\Bbb R^2$.) Therefore, the
points $\bar q_j(\lambda)=\bar q_j(0)$ ($j\in H_1$, $\lambda\in\Bbb R$)
are enclosed in bounded clusters, for they do not belong to $U$, see 
Lemma 5.3. Recall that, as the number $t_0$ varies, the whole set $U$ and
all of its connected components $U_0$ are moving in $\Bbb R^2$ at the
velocity $V_2-V_1$ (as the derivative with respect to $t_0$ shows), see
the term containing $t_0$ in (5.9). However, for the representatives
$\bar q_j+l$ ($l\in\Bbb Z^2$, $j\in H_1$) of the $H_1$-dynamics
$\{S_1^t\}$ (now the time parameter is $t_0$) it is impossible to remain
in a uniformly moving, bounded enclosure by Sub-lemma 5.5, for in that case all
velocities $v_j^\infty$ ($j\in H_1$) would be the same, contradicting to
the fact $|H_1|\ge2$ and the connectedness of the collision graph of $H_1$.
This proves Sub-lemma 5.14. \qed

\medskip

The joint conclusion of sub-lemmas 5.14 and 5.5 is that all velocities of the
internal flow $\{S_2^t\}$ of $H_2$ are parallel to the vector of periodicity
$l_0$ of the connected component $U_0$. Moreover, we constructed the lifting
$\bar q_i(\lambda)\in\Bbb R^2$ in such a way that the $l_0$-periodic strip
$U_0$ --- forbidden zone for the points $\bar q_i(\lambda)+\Bbb Z^2$
($i\in H_1$) --- moves at the velocity $V_2-V_1$, see the term containing
$t_0$ in (5.9). Since the lifting of the $H_1$-subsystem has no drift
(the ``baricenter'' is not moving when $t_0$ is changing, see the definition
of $\bar q_1(\lambda)=\bar q_1$ above),
we get that the relative velocity $V_2-V_1$
must also be parallel to $l_0$. This also means that the $\Bbb R^2$-lifting
of the internal flow $\{S_1^t\}$ is confined to an $l_0$-periodic, infinite
strip bounded by two translated copies of $U_0$. By using Sub-lemma 5.5 again,
we obtain that all velocities of the internal flow $\{S_1^t\}$ are also
parallel to the vector of periodicity $l_0$. \qed

\medskip

\subheading{Remark 5.14/a} Since 
$n(\lambda)=(t_0-\alpha\lambda)v_\infty+\lambda w_\infty$, the drift
(i. e. the average derivative of the positions with respect to the variable
$\lambda$) 
$$
(M_2)^{-1}\cdot\sum_{j\in H_2}m_j\frac{d}{d\lambda}\bar q_j(\lambda)
$$
of the subsystem $H_2$ is equal to $(W_2-W_1)-\alpha(V_2-V_1)$ where, as we 
recall, 
$$
V_i=(M_i)^{-1}\cdot\sum_{j\in H_i}m_jv_j^\infty, \quad
W_i=(M_i)^{-1}\cdot\sum_{j\in H_i}m_jw_j^\infty.
$$
Obviously, this drift must be parallel to the vector of periodicity $l_0$.
Since $V_2-V_1$ is parallel to $l_0$, we conclude that $W_2-W_1$ is also
parallel to $l_0$. This remark will be used later in this section. \endrem

\medskip

The second major case in ($3^o$) is

\subheading{Case (B): $\alpha=\beta$} Let us consider now the modified neutral
vector $w_\infty-\alpha v_\infty\in\Cal N_0(S^{\Bbb R}x_\infty)$. The advance
of both subsystems $H_1$ and $H_2$ is zero with respect to the neutral vector
$w_\infty-\alpha v_\infty$, thus
$$
\aligned
w_j^\infty-\alpha v_j^\infty=h_1, \quad \forall\, j\in H_1, \\
w_j^\infty-\alpha v_j^\infty=h_2, \quad \forall\, j\in H_2,
\endaligned
\tag 5.15
$$
for some vectors $h_1,\, h_2\in\Bbb R^2$. In other words, the effect of the
neutral translation by the vector $w_\infty-\alpha v_\infty$ on the
non-interacting groups $H_1$ and $H_2$ is that $H_i$ gets displaced 
(translated) by the vector $h_i$, $i=1,2$. Now there are again two sub-cases:

\medskip

\subheading{Sub-case B/1: $h_1\ne h_2$} In this case the result of the neutral
translation by the vector $n(\lambda)=\lambda(w_\infty-\alpha v_\infty)$
($\lambda\in\Bbb R$ is now the varying parameter) is that the relative
translation of the $H_2$ subsystem with respect to $H_1$ is
$\lambda(h_2-h_1)$ with the velocity $0\ne h_2-h_1\in\Bbb R^2$. The point is
that Lemma 5.3 is again readily applicable (so that $n(\mu,\lambda)$ is 
replaced by $n(\lambda)$), meaning that on the orbit 
$S^{\Bbb R}T_{n(\lambda),0}x_\infty$ of $T_{n(\lambda),0}x_\infty$ no proper
collision takes place between the groups $H_1$ and $H_2$. This fact has the
following consequence on the $\Bbb Z^2$-periodic, $\Bbb R^2$-lifting

$$
\left\{\bar q_i(t)+l\in\Bbb R^2\colon\; i\in H_1\cup H_2,\; l\in\Bbb Z^2,\;
t\in\Bbb R\right\}
\tag 5.16
$$
of the subsystem $H_1\cup H_2$ with the baricenter normalization
$\sum_{i\in H_1}m_i\dfrac{d}{dt}\bar q_i(t)=0$:

$$
\text{dist}\left(\bar q_i(t),\, \bar q_j(t)+\lambda(h_2-h_1)+l\right)
\ge 2r,
\tag 5.17
$$
for $i\in H_1$, $j\in H_2$, $t,\lambda\in\Bbb R$, $l\in\Bbb Z^2$. In other
words, the $2r$-wide, infinite strips with the direction of $h_2-h_1$
containing $\bar q_i(t)$ on their medium line ($i\in H_1$) are disjoint from
the similarly constructed infinite strips containing $\bar q_j(t)$ 
($j\in H_2$) on their medium line. Similarly to the closing part of the
discussion of Case (A), we conclude, first of all, that the relative
motion (drift) $V_2-V_1$ between $H_2$ and $H_1$ must be parallel to
$h_2-h_1$ and then, according to Sub-lemma 5.5, all velocities of the internal
dynamics $\{S_1^t\}$ and $\{S_2^t\}$ must also be parallel to $h_2-h_1$.
Since the $(h_2-h_1)$-parallel strips of width $2r$ are disjoint modulo
$\Bbb Z^2$, we immediately get that $h_2-h_1$ has a lattice direction,
and the shortest nonzero lattice vector $l_0$ parallel to $h_2-h_1$
has length at most $1/(4r)$. \qed

\medskip

\subheading{Remark 5.17/a} Let us observe that everything that has been said 
about the pair $(H_1,\, H_2)$ in Case B/1 can be repeated almost word-by-word
if one of the groups $H_i$, say $H_2$, has only one element. This remark
will have a particular relevance later in this section. \endrem

\medskip

\subheading{Sub-case B/2: $h_1=h_2$} In this situation the united subsystem
$H_1\cup H_2$ gets uniformly translated by the vector $h_1=h_2$ under the
action of the neutral vector $w_\infty-\alpha v_\infty$. This is an open 
possibility, indeed, and nothing else can be said about it.

\medskip

Now we are in the position of quickly finishing the discussion of ($3^o$).
Recall that $s=k$ ($\ge2$), i. e. $|H_i|\ge2$ for $i=1,\dots,k$. Consider
the advances $\alpha_i=\alpha(H_i)$ of the subsystems $H_1,\dots,H_k$
with respect to the limiting neutral vector 
$w_\infty=\lim_{n\to\infty}w_{t_n}\in\Cal N_0(S^{\Bbb R}x_\infty)$.
Unfortunately, we again have to distinguish between two cases. 

\medskip

\subheading{Case I. Not all $\alpha_i$'s are the same, e. g.
$\alpha_1\ne\alpha_2$} In this situation the result of Case (A) above says
that $V_1-V_2$ and all velocities of the internal flows $\{S_1^t\}$ and
$\{S_2^t\}$ are parallel to the same nonzero lattice vector $l_0\in\Bbb Z^2$.
For any other subsystem $H_i$ ($i>2$) we have that $\alpha_i\ne\alpha_1$ or
$\alpha_i\ne\alpha_2$. Assume that $\alpha_i\ne\alpha_1$. The result of 
Case (A) above says that $V_i-V_1$ and all velocities of the internal flows 
$\{S_1^t\}$ and $\{S_i^t\}$ are parallel to the same nonzero lattice vector 
$l_1\in\Bbb Z^2$. The common presence of the flow $\{S_1^t\}$ in these
statements shows that $l_0=l_1$ (or, at least they are parallel to each
other). Summarizing these results, we finish the discussion of Case I by
concluding that all average velocities $V_i$ and all velocities of the
internal flows $\{S_i^t\}$ ($i=1,\dots,k$) are parallel to the same (nonzero)
lattice vector $l_0$ whose magnitude is at most $1/(4r)$. \qed

\medskip

\subheading{Case II. $\alpha_1=\alpha_2=\dots=\alpha_k=:\alpha$} Consider,
as in Case (B) above, the neutral vector
$w_\infty-\alpha v_\infty\in\Cal N_0(S^{\Bbb R}x_\infty)$. The advance of
each $H_i$ with respect to $w_\infty-\alpha v_\infty$ is zero, so (5.15)
applies:
$$
w_j^\infty-\alpha v_j^\infty=h_i \text{ for } j\in H_i,\; i=1,2,\dots,k.
\tag 5.18
$$
Since $w_\infty-\alpha v_\infty\ne 0$ and $\sum_{i=1}^k M_ih_i=0$
($M_i=\sum_{j\in H_i}m_j$), we conclude that not all vectors
$h_1,\dots,h_k\in\Bbb R^2$ are the same, e. g. $h_1\ne h_2$. Then for every
$i>2$ we have either $h_i\ne h_1$, or $h_i\ne h_2$, say $h_i\ne h_1$, and the
result of Case B/1 above applies to the pairs of subsystems $(H_1,H_2)$ and
$(H_1,H_i)$. Quite similarly to the discussion of Case I above (but referring
in it to Case B/1, instead of Case (A)) we get that all average velocities
$V_1,\dots,V_k$ and all velocities of the internal flows $\{S_i^t\}$
($i=1,\dots,k$) are parallel to the vector $h_2-h_1\ne0$. Recall that 
$h_2-h_1$ has a lattice direction, and the shortest (nonzero) lattice vector
$l_0$ parallel to $h_2-h_1$ has magnitude at most $1/(4r)$, thus completing
the proof of Key Lemma 5.4 in the case ($3^o$). \qed

\medskip

\subheading{$4^o$ The general case $1\le s<k$} In this case the subsystems 
$H_i$ with $|H_i|\ge2$ (i. e. $i\le s$) coexist with the subsystems $H_i$
for which $|H_i|=1$ ($i>s$). In order to simplify the notations we assume
that $H_i=\{i\}$ for $i=s+1,s+2,\dots,k$.

Consider the advances $\alpha_j=\alpha(H_j)$ of the subsystems $H_j$,
$1\le j\le s$, with respect to the limiting neutral vector $w_\infty$.
Unfortunately, we again have to distinguish between three major situations.

\medskip

\subheading{Case I. Not all of $\alpha_1,\dots,\alpha_s$ are equal,
say, $\alpha_1\ne\alpha_2$} Let us observe, first of all, that the whole
machinery of ($3^o$) applies to the united subsystem
$H_1\cup H_2\cup\dots\cup H_s$, showing that there exists a nonzero lattice 
vector $l_0\in\Bbb Z^2$ so that all velocities of the internal dynamics
$\{S_i^t\}$ ($1\le i\le s$) and all relative velocities of the baricenters
$V_i-V_j$ ($1\le i,j\le s$) are parallel to $l_0$. (The second part of this
statement is clearly equivalent to saying that the average velocities
$V'_i$ of the subsystems $H_i$ ($1\le i\le s$) are parallel to $l_0$, 
provided that these average velocities are observed from a reference system
attached to the baricenter of $H_1\cup H_2\cup\dots\cup H_s$.) 

Let us turn our attention to a one-disk subsystem $H_i=\{i\}$, 
$s+1\le i\le k$, $i$ is fixed.
Just like in (5.18), the advance of the subsystem $H_j$
($j\le s$) with respect to the neutral vector $w_\infty-\alpha_j v_\infty$
is zero, therefore the whole subsystem $H_j$ gets translated by the same
vector $h_j\in\Bbb R^2$ under the action of $w_\infty-\alpha_j v_\infty$:
$$
w_l^\infty-\alpha_j v_l^\infty=h_j, \quad l\in H_j,\; j=1,\dots,s.
\tag 5.19
$$
Consider now the vectors of displacement
$w_i^\infty-\alpha_j v_i^\infty=h'_j\in\Bbb R^2$, $j=1,\dots,s$, 
the index $i$ is fixed,
$s+1\le i\le k$. If $h_j-h'_j\ne0$ for at least one $j\le s$, then the result
of Case B/1 of ($3^o$) applies to the pair of subsystems $(H_j,\, H_i)$
(see Remark 5.17/a), thus we have that the relative velocity $V_j-V_i$
of the baricenters is parallel to the fixed lattice vector $l_0\in\Bbb Z^2$.
This is the most we can  prove for the motion of $H_i$ relative to the motion
of $H_1\cup H_2\cup\dots\cup H_s$, for if we had such a result for every
$i$ ($s+1\le i\le k$), then the statement of the key lemma would follow.

The unpleasant situation with $H_i$ is when
$$
w_i^\infty-\alpha_j v_i^\infty=h_j=w_l^\infty-\alpha_j v_l^\infty,\;
l\in H_j,\; j=1,\dots,s,
\tag 5.20
$$
$s+1\le i\le k$, $i$ is fixed.

With $i$ and $j$ fixed, let us average (5.20) with respect to the weights
$m_l$ ($l\in H_j$) of the subsystem $H_j$. We obtain
$$
W_j=W_i+\alpha_j(V_j-V_i), \quad j=1,\dots,s.
\tag 5.21
$$
$s+1\le i\le k$, $i$ is fixed. Recall that
$V_i=v_i^\infty$ and $W_i=w_i^\infty$ for the one-disk subsystem $H_i=\{i\}$.

We again have to distinguish between two sub-cases.

\medskip

\subheading{Case I/a. Not all average velocities $V_1,\dots,V_s$ are the same}

\medskip

\proclaim{Sub-lemma 5.22} There is a pair of indices $1\le j_1,\, j_2\le s$
for which $\alpha_{j_1}\ne\alpha_{j_2}$ and $V_{j_1}\ne V_{j_2}$.
\endproclaim

\subheading{Proof} As a matter of fact, this sub-lemma is trivial. Indeed,
if $V_{j_1}$ were equal to $V_{j_2}$ whenever $\alpha_{j_1}\ne\alpha_{j_2}$
($1\le j_1,\, j_2\le s$), then we would have, first of all, $V_1=V_2$,
since $\alpha_1\ne\alpha_2$ by the assumption of Case I. Secondly, for every
$j=3,4,\dots,s$ either $\alpha_j\ne\alpha_1$ or $\alpha_j\ne\alpha_2$,
thus proving $V_j=V_1=V_2$ for $j=3,4,\dots,s$, contradicting to the 
assumption of I/a. \qed

\medskip

By taking the difference of (5.21) for $j_1$ and $j_2$, and also using 
Remark 5.14 for the pair of subsystems $(H_{j_1},\, H_{j_2})$, we get
$$
\aligned
c(V_{j_1}-V_{j_2})=W_{j_1}-W_{j_2}=\alpha_{j_1}(V_{j_1}-V_{i})-
\alpha_{j_2}(V_{j_2}-V_{i}) \\
=\alpha_{j_1}(V_{j_1}-V_{j_2})+(\alpha_{j_1}-\alpha_{j_2})(V_{j_2}-V_{i}),
\endaligned
\tag 5.23
$$
for some scalar $c$. Since $V_{j_1}-V_{j_2}\parallel l_0$ and
$\alpha_{j_1}-\alpha_{j_2}\ne0$, we obtain that $V_{j_2}-V_{i}$ is also
parallel to the lattice vector $l_0$, precisely what we wanted to prove in
Case I.

\medskip

\subheading{Case I/b. $V_1=V_2=\dots=V_s=:V$} Now formula (5.21) says that
$$
W_j=W_i+\alpha_j(V-V_i), \quad j=1,\dots,s,
\tag 5.24
$$
the index $i$ is fixed, $s+1\le i\le k$. Take the difference of (5.24) 
for $j=1$ and $j=2$:
$$
W_1-W_2=(\alpha_1-\alpha_2)\cdot(V-V_i).
\tag 5.25
$$
Recall that $\alpha_1\ne\alpha_2$ and, by Remark 5.14/a, $W_1-W_2$ is parallel
to $l_0$. Therefore, the relative velocity $V_i-V$ also proves to be parallel
to the lattice vector $l_0$, the result we just wanted to prove for the
subsystem $H_i$ in Case I. Thus Key Lemma 5.4 has been proved in Case I of
($4^o$).

\medskip

\subheading{Case II. $\alpha_1=\alpha_2=\dots=\alpha_s=:\alpha$, but not all
vectors $h_1,\dots,h_s$ in (5.19) are the same} Assume that $h_1\ne h_2$.
Then the method of Sub-case B/1 of ($3^o$) applies to the subsystem 
$H_1\cup\dots\cup H_s$, showing again that there exists a nonzero lattice
vector $l_0$ such that all velocities of the flows $\{S_i^t\}$
($1\le i\le s$) and all relative velocities $V_{j_1}-V_{j_2}$ 
($1\le j_1,\, j_2\le s$)
are parallel to $l_0$. Just as in Case I above, consider again a one-disk
subsystem $H_i=\{i\}$, $s+1\le i\le k$. We need to show that some (or any)
of the velocities $V_i-V_j$ ($j=1,\dots,s$) is parallel to $l_0$. Consider
the vector $h_i=w_i^\infty-\alpha v_i^\infty\in\Bbb R^2$. This vector should
differ from $h_1$ or $h_2$. Assume that $h_i\ne h_1$. In this situation the
method and result of Sub-case B/1 of ($3^o$) again applies to the pair of
subsystems $(H_1,\, H_i)$ (see Remark 5.17/a), and we obtain that $V_i-V_1$
is parallel to $l_0$. This step finishes the proof of Key Lemma 5.4 in 
Case II of ($4^o$). \qed

\medskip

\subheading{Case III. $\alpha_1=\alpha_2=\dots=\alpha_s=:\alpha$, and
$h_1=h_2=\dots=h_s=:h$ in (5.19)} We can assume that 
$h_1=\dots=h_t=h$ and $h_i\ne h$ for $t+1\le i\le k$. Due to the relation
$\sum_{i=1}^k M_ih_i=0$ (and to the fact that $w_\infty\ne\alpha v_\infty$),
we have that $s\le t\le k-1$. Select and fix an
arbitrary index $i\in\{t+1,\dots,k\}$, and study the relative motion of the
subsystems $H^*=:H_1\cup H_2\cup\dots\cup H_t$ and $H_i=\{i\}$. Since the
neutral vector $w_\infty-\alpha v_\infty$ translates the whole subsystem 
$H^*$ by the same vector $h$ and it translates the one-disk subsystem
$H_i=\{i\}$ by a different vector $h_i$, the method and result of Sub-case B/1
of ($3^o$) again applies to the pair $(H^*,\, H_i)$, and we obtain that
there exists a nonzero lattice vector $l_0\in\Bbb Z^2$ so that all velocities
of the internal dynamics $\{S_j^t\}$ ($j=1,\dots,s$), all relative velocities
$V_{j_1}-V_{j_2}$ ($j_1,\, j_2\in\{1,\dots,t;i\}$), and $h_i-h$ are parallel 
to $l_0$. Due to the common presence of the internal dynamics $\{S_1^t\}$,
the same thing can be said about any other index $i\in\{t+1,\dots,k\}$
with the same direction vector $l_0$. This finishes the proof of Key Lemma 5.4
in the last remaining case of ($4^o$), thus completing the proof of 5.4.
We note that every nonzero lattice vector $l_0\in\Bbb Z^2$ that emerged in
this proof had the property $||l_0||\le 1/(4r)$, thus ensuring the finiteness
of the family $\{l_0,l_1,\dots,l_p\}$ in Key Lemma 5.4. \qed

\bigskip \bigskip

\heading
\S6. Non-Existence of Separating Manifolds \\
Part C: Topological Arguments
\endheading

\bigskip \bigskip

Given any nonzero vector $l_0\in\Bbb Z^2$ ($||l_0||\le \dfrac{1}{4r}$,
as always), consider the one-dimensional sub-torus
$T(l_0)=\left\{\lambda l_0|\; \lambda\in\Bbb R\right\}/\Bbb Z^2$ of
$\Bbb T^2=\Bbb R^2/\Bbb Z^2$, and define the following subset $L(l_0)$
of the phase space $\bold M$:

$$
L(l_0)=\left\{x\in\bold M|\; v_i(S^tx)\parallel l_0\;\; \forall t\in\Bbb R,
\;\; i=1,\dots,N\right\}.
\tag 6.1
$$
The set $L(l_0)$ is obviously flow-invariant. We will see below that $L(l_0)$
is a compact subset of the (compact) phase space $\bold M$. As a consequence
of Key Lemma 5.4 and Remark 5.4/a, we have that for any separating manifold
$J\subset\bold M\setminus\partial\bold M$ (with the dimension-minimizing 
property, see properties (0)--(3) at the end of \S3) and for every phase
point $x_0\in J$ --- with a non-singular forward orbit 
$S^{[0,\infty)}x_0$ --- there exists a lattice vector $l_0\in\Bbb Z^2$
($0<||l_0||\le\dfrac{1}{4r}$) such that the $\Omega$-limit set 
$\Omega(x_0)$ of $x_0$ is contained in $L(l_0)$.

Let us briefly describe first the structure of the set $L(l_0)$. For any point
$x\in L(l_0)$ we define the partition $\Cal P=\Cal
P(x)=\left\{H_1,\dots,H_k\right\}=\left\{H_1(x),\dots, H_k(x)\right\}$
($k=k(x)$) of the full vertex set $\{1,2,\dots,N\}$ into the connected
components of the collision graph $\Cal G(S^{\Bbb R}x)$ of the orbit
$S^{\Bbb R}x$, just as we did in 5.2 above. (Recall that in the construction
of the collision graph $\Cal G(S^{\Bbb R}x)$ we only consider the proper,
i. e. non-tangential collisions.) For definiteness, let the labeling of the
sets $H_i=H_i(x)$ ($i=1,\dots,k(x)$) follow the pattern that for $i<j$
the smallest element of $H_i$ precedes the smallest element of $H_j$.
Denote the open, tubular $r$-neighborhood (in $\Bbb T^2$) of the coset
$T(l_0)+q_i(x)$ by $T_i=T_i(x)$, $x\in L(l_0)$, $i=1,\dots,N$. Thanks to the
invariance of the relation $v_i(S^tx)\parallel l_0$ ($\forall t\in\Bbb R$),
we obtain that $T_i(S^tx)=T_i(x)$ for all $t\in\Bbb R$.

A simple observation of the trajectory of a phase point $x\in L(l_0)$
reveals the following facts about any pair of {\it distinct} indices
$i_1,\, i_2\in\{1,2,\dots,N\}$:

$$
T_{i_1}(x)=T_{i_2}(x),\text{ if } i_1,\, i_2\in H_r(x)\text{ for some }
1\le r\le k(x),
\tag 6.2
$$

$$
\aligned
&T_{i_1}(x)\cap T_{i_2}(x)=\emptyset \\
&\text{if } i_1\in H_{r_1}(x),\;
i_2\in H_{r_2}(x),\; r_1\ne r_2,\text{ and }\max\left\{|H_{r_1}(x)|,\,
|H_{r_2}(x)|\right\}\ge2,
\endaligned
\tag 6.3
$$

$$
T_{i_1}(x)\cap T_{i_2}(x)=\emptyset\text{ or } v_{i_1}(x)=v_{i_2}(x)
\text{ if } \max\left\{|H_{r_1}(x)|,\, |H_{r_2}(x)|\right\}=1.
\tag 6.4
$$
It is clear from these formulas that the set of phase points $x\in L(l_0)$
with a given partition $\Cal P=\Cal P(x)$ is closed (i. e. compact) and,
henceforth, the entire set $L(l_0)$ is also compact.

The following proposition is, in fact, a simple consequence of the Baire
category theorem.

\medskip

\proclaim{\bf Proposition 6.5} Assume that
$J\subset\bold M\setminus\partial\bold M$ is a separating
manifold with the dimension minimalizing property, i. e. properties
(0)---(3) from the end of \S3 hold true. We claim that there exists
a lattice vector $l_0\in\Bbb Z^2$ (with $0<||l_0||\le 1/4r$) such that for
every $\epsilon_0>0$ there is an open subset $\emptyset\ne G\subset J$
and a threshold $t_0>0$ for which

$$
\aligned
d\left(S^ty,\, L(l_0)\right)&\le\epsilon_0, \\
\lim_{\tau\to\infty}d\left(S^{\tau}y,\, L(l_0)\right)&=0
\endaligned
\tag 6.6
$$
for all $y\in G$ with a non-singular forward orbit $S^{(0,\infty)}y$
and for all $t\ge t_0$.
\endproclaim

\medskip

\subheading{\bf Proof}
Denote by $\Cal S_J=\Cal S^+\cap J$ the set of all phase
points $y\in J$ with a singular forward orbit $S^{(0,\infty)}y$. It follows
from property (3) at the end of \S3 that the set $\Cal S_J$ is an 
$F_\sigma$ set (i. e. a countable union of closed sets) of zero measure in $J$.
For fixed $\epsilon_0>0$, $t_0>0$, and $l_0\in\Bbb Z^2\setminus\{0\}$ the set

$$
F(\epsilon_0,t_0,l_0)=\left\{y\in J\setminus\Cal S_J|\; \forall t\ge t_0
\quad d(S^ty,\, L(l_0))\le\epsilon_0\right\}
$$
is closed in the dense, $G_\delta$ subset (countable intersection of open
sets) $J\setminus\Cal S_J$ of $J$. It follows from the results of earlier
sections (Key Lemma 5.4, Remark 5.4/a) that

$$
J\setminus\Cal S_J=\bigcup\Sb t_0>0 \\ l_0\in\Bbb Z^2\setminus\{0\} \\
||l_0||\le 1/4r \endSb F(\epsilon_0,t_0,l_0)
\tag 6.7
$$
for any fixed $\epsilon_0>0$. The space $J\setminus\Cal S_J$ is a $G_\delta$
subspace of $J$ and, being such, it is completely metrizable, see Theorem
4.3.23 in [E(1977)]. Consequently, in the topological space 
$J\setminus\Cal S_J$ the Baire category theorem is applicable. The union in
(6.7) is monotonic in $t_0$, thus --- by the Baire theorem --- there exist a
$t_0>0$ and a direction vector $l_0\in\Bbb Z^2\setminus\{0\}$
($||l_0||\le 1/4r$) such that the set $F(\epsilon_0,t_0,l_0)$ has a non-empty
interior in $J\setminus\Cal S_J$, i. e. there is an open set
$\emptyset\ne G\subset J$ such that

$$
\aligned
d\left(S^ty,\, L(l_0)\right)&\le\epsilon_0, \\
\lim_{\tau\to\infty}d\left(S^{\tau}y,\,\bigcup\Sb l\in\Bbb Z^2 \\
0<||l||\le 1/4r \endSb L(l)\right)&=0
\endaligned
$$
for all $y\in J\setminus\Cal S_J$ and $t\ge t_0$. However, the compact
components $L(l)$ of

$$
\bigcup\Sb l\in\Bbb Z^2 \\ 0<||l||\le 1/4r \endSb L(l)
$$
are mutually disjoint (Needless to say, we only consider mutually
non-parallel vectors $l$!), so if the given number $\epsilon_0>0$ is
selected to be sufficiently small, then --- by the already proved first
inequality of (6.6) --- the orbit $S^ty$ ($y\in G\setminus\Cal S_J$) cannot
converge to any component $L(l)$ other than $L(l_0)$.

Finally, for all smaller values $\epsilon'_0<\epsilon_0$ we can repeat the
above argument by restricting ourselves to the open set $G$ instead of the
entire $J$ and, in that way, the direction vector $l$ of the limiting set
$L(l)$ will remain the above vector $l_0$ all the time (i. e. for every
$\epsilon'_0$). This finishes the proof of the proposition. \qed

\medskip

\subheading{\bf Corollary 6.8} By replacing the manifold $J$ by $S^{t_0}(G)$,
we can assume that the statement of Proposition 6.5 holds true for the entire
separating manifold $J$ with the threshold $t_0=0$. What is even more, the
condition on the non-singularity of $S^{(0,\infty)}y$ ($y\in J$) can be
dropped, as the following argument shows:

The forward orbit of a phase point $y\in\Cal S_J$ has several branches, see
\S2.3 above, or \S2 of [Sim(1992-A)]. However, each of these branches is
actually the limit of forward orbits of phase points 
$y_n\in J\setminus\Cal S_J$ ($n\to\infty$). Consequently, the first line of
(6.6) (the inequality) readily generalizes to the singular phase points
$y\in G\cap\Cal S_J$. Even if it might seem appealing, throughout the entire
proof of the Theorem we will not need this additional result about singular
phase points. We will be exclusively dealing with phase points 
$y\in\bold M\setminus\partial\bold M$ with a nonsingular forward orbit
$S^{(0,\infty)}y$. \qed

In the upcoming two sections we are going to prove that --- contrary
to the statement of Proposition 6.5 --- the compact set $L(l_0)$ 
($l_0\in\Bbb Z^2$ fixed, $0<||l_0||\le 1/4r$) cannot attract any separating
manifold $J\subset\bold M\setminus\partial\bold M$. The accomplishment of such
a proof will complete the proof of the Theorem of this paper.

\bigskip \bigskip

\heading
\S7. Non-Existence of Separating Manifolds \\
Part D: Transversality
\endheading

\bigskip \bigskip

Given a codimension-one, locally
flow-invariant, smooth sub-manifold $J\subset\bold M$,
consider a normal vector $n_0=(z,w)$ ($\ne 0$) of $J$ at the phase point
$y\in J$, i. e. for any tangent vector 
$(\delta q,\, \delta v)\in\Cal T_y\bold M$ the relation 
$(\delta q,\, \delta v)\in\Cal T_yJ$ is true if and only if 
$\langle\delta q,z\rangle+\langle\delta v,w\rangle=0$. Here 
$\langle\, .\, ,\, .\, \rangle$ is the scalar product corresponding to the
mass metric, that is,
$\langle a,\, b\rangle=\sum_{i=1}^{N}m_i\langle a_i,\, b_i\rangle$. Let us
determine first the time-evolution $n_0\longmapsto n_t$ ($t>0$) of this normal
vector as time $t$ elapses. If there is no collision on the orbit segment
$S^{[0,t]}y$, then the relationship between
$(\delta q,\, \delta v)\in\Cal T_y\bold M$ and 
$(\delta q',\, \delta v')=\left(DS^t\right)(\delta q,\, \delta v)$ is
obviously

$$
\aligned
\delta v'&=\delta v, \\
\delta q'&=\delta q+t\delta v,
\endaligned
\tag 7.1
$$
from which we obtain that

$$
\aligned
(\delta q',\, \delta v')\in\Cal T_{y'}J&\Leftrightarrow\langle\delta q'-t
\delta v',\, z\rangle+\langle\delta v',\, w\rangle=0 \\
&\Leftrightarrow
\langle\delta q',\, z\rangle+\langle\delta v',\, w-tz\rangle=0.
\endaligned
$$
This means that $n_t=(z,\, w-tz)$. It is always very useful to consider
the quadratic form $Q(n)=Q((z,w))=:\langle z,w\rangle$ associated with the
normal vector $n=(z,w)\in\Cal T_y\bold M$ of $J$ at $y$. $Q(n)$ is the
so called ``infinitesimal Lyapunov function'', see [K-B(1994)] or part
A.4 of the Appendix in [Ch(1994)]. For a detailed exposition of the 
relationship between the quadratic form $Q$, the relevant symplectic geometry
and the dynamics, please see [L-W(1995)].

\medskip

\subheading{\bf Remark} Since the normal vector $n=(z,w)$ of $J$ is only
determined up to a nonzero scalar multiplier, the value $Q(n)$ is only
determined up to a positive multiplier. However, this means that the sign
of $Q(n)$ (which is the utmost important thing for us) is uniquely
determined. This remark will gain a particular importance in the near
future.

\medskip

\noindent
From the above calculations we get that 

$$
Q(n_t)=Q(n_0)-t||z||^2\le Q(n_0).
\tag 7.2
$$

The next question is how the normal vector $n$ of $J$ gets transformed
$n^-\mapsto n^+$ through a collision (reflection) at time $t=0$? Elementary
geometric considerations show (see Lemma 2 of [Sin(1979)], or formula
(2) in \S3 of [S-Ch(1987)]) that the linearization of the flow

$$
\left(DS^t\right)\Big|_{t=0}:\; (\delta q^-,\, \delta v^-)\longmapsto
(\delta q^+,\, \delta v^+)
$$
is given by the formulas

$$
\aligned
\delta q^+&=R\delta q^-, \\
\delta v^+&=R\delta v^-+2\cos\phi RV^*KV\delta q^-,
\endaligned
\tag 7.3
$$
where the operator $R:\; \Cal T_q\bold Q\to \Cal T_q\bold Q$ is the orthogonal
reflection (with respect to the mass metric) across the tangent hyperplane
$\Cal T_q\partial\bold Q$ of $\partial\bold Q$ at $q\in \partial\bold Q$
($y^-=(q,v^-)\in\partial\bold M$, $y^+=(q,v^+)\in\partial\bold M$), 
$V:\; (v^-)^\perp\to\Cal T_q\partial\bold Q$ is the $v^-$-parallel projection
of the ortho-complement hyperplane $(v^-)^\perp$ onto
$\Cal T_q\partial\bold Q$, $V^*:\; \Cal T_q\partial\bold Q\to (v^-)^\perp$
is the adjoint of $V$, i. e. it is the projection of $\Cal T_q\partial\bold Q$
onto $(v^-)^\perp$ being parallel to the normal vector $\nu(q)$ of
$\partial\bold Q$ at $q\in\partial\bold Q$, 
$K:\; \Cal T_q\partial\bold Q\to \Cal T_q\partial\bold Q$ is the second
fundamental form of $\partial\bold Q$ at $q$ and, finally, 
$\cos\phi=\langle\nu(q),\, v^+\rangle$ is the cosine of the angle $\phi$
subtended by $v^+$ and the normal vector $\nu(q)$. For the formula (7.3),
please also see the last displayed formula of \S1 in [S-Ch(1982)], or
(i) and (ii) of Proposition 2.3 in [K-S-Sz(1990)]. We note that it is enough
to deal with the tangent vectors 
$(\delta q^-,\, \delta v^-)\in(v^-)^\perp\times(v^-)^\perp$
($(\delta q^+,\, \delta v^+)\in(v^+)^\perp\times(v^+)^\perp$), for the
manifold $J$ under investigation is supposed to be flow-invariant, so any
vector $(\delta q,\, \delta v)=(\alpha v,\, 0)$ ($\alpha\in\Bbb R$) is
automatically inside $\Cal T_yJ$. The backward version (inverse) 

$$
\left(DS^t\right)\Big|_{t=0}:\; (\delta q^+,\, \delta v^+)\mapsto
(\delta q^-,\, \delta v^-)
$$
can be deduced easily from (7.3):

$$
\aligned
\delta q^-&=R\delta q^+, \\
\delta v^-&=R\delta v^+-2\cos\phi RV_1^*KV_1\delta q^+,
\endaligned
\tag 7.4
$$
where $V_1:\; (v^+)^\perp\to\Cal T_q\partial\bold Q$ is the $v^+$-parallel
projection of $(v^+)^\perp$ onto $\Cal T_q\partial\bold Q$. By using formula
(7.4), one easily computes the time-evolution $n^-\longmapsto n^+$
of a normal vector $n^-=(z,w)\in\Cal T_{y^-}\bold M$ of $J$ if a collision
$y^-\longmapsto y^+$ takes place at time $t=0$:

$$
\aligned
(\delta q^+,\, \delta v^+)\in\Cal T_{y^+}J\Leftrightarrow\langle R\delta q^+,
\, z\rangle+\langle R\delta v^+-2\cos\phi RV_1^*KV_1\delta q^+,\,
w\rangle &=0 \\
\Leftrightarrow\langle\delta q^+,\, Rz-2\cos\phi V_1^*KV_1Rw\rangle+
\langle\delta v^+,\, Rw\rangle &=0.
\endaligned
$$
This means that

$$
n^+=\left(Rz-2\cos\phi V_1^*KV_1Rw,\, Rw\right)
\tag 7.5
$$
if $n^-=(z,\, w)$. It follows that

$$
\aligned
Q(n^+)&=Q(n^-)-2\cos\phi\langle V_1^*KV_1Rw,\, Rw\rangle \\
&=Q(n^-)-2\cos\phi\langle KV_1Rw,\, V_1Rw\rangle\le Q(n^-).
\endaligned
\tag 7.6
$$
Here we used the fact that the second fundamental form $K$ of 
$\partial\bold Q$ at $q$ is positive semi-definite, which just means that the
billiard system is semi-dispersive.

The last simple observation on the quadratic form $Q(n)$
regards the involution $I:\; \bold M\to\bold M$, $I(q,v)=(q,-v)$ 
corresponding to the time reversal. If $n=(z,w)$ is a normal vector of $J$
at $y$, then, obviously, $I(n)=(z,-w)$ is a normal vector of $I(J)$ at
$I(y)$ and

$$
Q\left(I(n)\right)=-Q(n).
\tag 7.7
$$

By switching --- if necessary --- from the separating manifold $J$ to 
$I(J)$, and by taking a suitable remote image $S^t(J)$ ($t>>1$), in the
spirit of (7.2), (7.6)--(7.7) we can assume that

$$
Q(n)\le c_0'<0
\tag 7.8
$$
uniformly for every {\it unit} normal vector $n\in\Cal T_y\bold M$ of
$J$ at any phase point $y\in J$.

\medskip

\subheading{\bf Remark 7.9} There could be, however, a little difficulty in
achieving the inequality $Q(n)<0$, i. e. (7.8). Namely, it may happen that
$Q(n_t)=0$ for every $t\in\Bbb R$. According to (7.2), the equation $Q(n_t)=0$
($\forall\, t\in\Bbb R$) implies that 
$n_t=:(z_t,\, w_t)=(0,\, w_t)$ for all $t\in\Bbb R$ and, moreover, in the view
of (7.5), $w_t^+=Rw_t^-$ is the transformation law at any collision
$y_t=(q_t,\, v_t)\in\partial\bold M$. Furthermore, at every collision
$y_t=(q_t,\, v_t)\in\partial\bold M$ the projected tangent vector
$V_1Rw_t^-=V_1w_t^+$ lies in the null space of the operator $K$ 
(see also (7.5)), and this means that $w_0$ is a neutral vector for the
entire trajectory $S^{\Bbb R}y$, i. e. $w_0\in\Cal N\left(S^{\Bbb R}y\right)$.
(For the notion of neutral vectors and $\Cal N\left(S^{\Bbb R}y\right)$,
cf. \S\S2.4 above.) On the other hand, this is impossible
for the following reason: Any tangent vector $(\delta q,\delta v)$ from the
space $\Cal N\left(S^{\Bbb R}y\right)\times\Cal N\left(S^{\Bbb R}y\right)$
is automatically tangent to the separating manifold $J$ (as a direct
inspection shows), thus for any normal vector $n=(z,w)\in\Cal T_y\bold M$
of a separating manifold $J$ one has

$$
(z,\, w)\in\Cal N\left(S^{\Bbb R}y\right)^\perp\times\Cal N\left(
S^{\Bbb R}y\right)^\perp.
\tag 7.10
$$
The membership in (7.10) is, however, impossible with a nonzero vector
$w\in\Cal N\left(S^{\Bbb R}y\right)$. \qed

\medskip

\subheading{\bf Singularities} 

\medskip

Consider a smooth, connected piece
$\Cal S\subset\bold M$ of a singularity manifold corresponding to a
singular (tangential or double) reflection {\it in the future}. Such a
manifold $\Cal S$ is locally
flow-invariant and has one codimension, so we can
speak about its normal vectors $n$ and the uniquely determined sign of
$Q(n)$ for $0\ne n\in\Cal T_y\bold M$, $y\in\Cal S$, $n\perp\Cal S$
(depending on the foot point, of course). Consider first a phase point
$y^-\in\partial\bold M$ right before the singular reflection that is
described by $\Cal S$. It follows from the proof of Lemma 4.1 of
[K-S-Sz(1990)] and Sub-lemma 4.4 therein that at 
$y^-=(q,\, v^-)\in\partial\bold M$ any tangent vector 
$(0,\, \delta v)\in\Cal T_{y^-}\bold M$ lies actually in 
$\Cal T_{y^-}\Cal S$ and, consequently, the normal vector
$n=(z,w)\in\Cal T_{y^-}\bold M$ of $\Cal S$ at $y^-$ necessarily has the
form $n=(z,0)$, i. e. $w=0$. Thus $Q(n)=0$ for any normal vector
$n\in\Cal T_{y^-}\bold M$ of $\Cal S$. According to the monotonicity 
inequalities (7.2) and (7.6) above,

$$
Q(n)>0
\tag 7.11
$$
for any phase point $y\in\Cal S$ of a future singularity manifold $\Cal S$.
As an immediate consequence of the inequalities (7.8) and (7.11), the summary
of this section is

\medskip

\proclaim{\bf Proposition 7.12} In some neighborhood of any phase point
$x_0\in J$ of a separating manifold $J$ (fulfilling (7.8) and conditions
(0)---(3) at the end of \S3 above) the manifold $J$ is uniformly 
transversal to any future singularity manifold $\Cal S$. Here the phrase
``uniform transversality'' means that in some open neighborhood $U_0$
of $x_0$ it is true that all possible angles 
$\alpha=\angle\left(\Cal T_y\Cal S,\Cal T_zJ\right)$ subtended by a tangent
space $\Cal T_y\Cal S$ of a future singularity ($y\in U_0\cap\Cal S$,
no matter what the order of the singularity) and a tangent space $\Cal T_zJ$
are separated from zero.
\endproclaim

\bigskip \bigskip

\heading
\S8. Non-Existence of Separating Manifolds \\
Part E: Dynamical-Geometric Considerations
\endheading

\bigskip \bigskip

\subheading{\bf The foliation} 

\medskip

By using propositions 6.5 and 7.12, for any fixed, small number $\epsilon_0>0$
let us consider a separating manifold $J\subset\bold M\setminus\partial\bold M$
enjoying all properties (0)---(3) from the end of \S3 so that also the 
transversality property (the statement of Proposition 7.12) holds true for $J$ 
and, finally,

$$
\aligned
d\left(S^tJ,\, L(l_0)\right)&\le\epsilon_0\;\; \forall t\ge 0, \\
\lim_{\tau\to\infty}d\left(S^{\tau}y,\, L(l_0)\right)&=0
\endaligned
\tag 8.1
$$
for all $y\in J\setminus\Cal S_J$. (Recall that $\Cal S_J$ denotes the set of
all phase points $y\in J$ with a singular forward orbit $S^{(0,\infty)}y$.)
The validity of the following proposition follows directly from 
Proposition 7.12 by also using the actual inequalities (7.8) and (7.11)
leading to 7.12.

\medskip

\proclaim{\bf Proposition 8.2} For any separating manifold $J$ (enjoying all 
properties described above) there exists a non-empty, open subset $G$ of $J$
that admits a smooth foliation $G=\bigcup_{i\in I}F_i$ by the curves $F_i$ 
with the following properties:

(1) The smooth curves $F_i$ are uniformly transversal to all future
singularities $\Cal S$, where uniformity is meant just as in Proposition 7.12;

(2) The curves $F_i$ are uniformly convex in the sense that for any (nonzero)
tangent vector $\tau=(\delta q,\, \delta v)\in\Cal T_yF_i$ ($y=(q,v)\in F_i$)
it is true that $\delta q\perp v$, $\delta v\perp v$, and
$\langle\delta q,\delta v\rangle/\Vert\delta q\Vert^2\ge c_0>0$
with some constant $c_0>0$ depending only on $G$;

(3) Write the components of the tangent vector
$0\ne\tau=(\delta q,\delta v)\in\Cal T_yF_i$ in the form 
$\delta q=\delta q^0+\delta q^\perp$, $\delta v=\delta v^0+\delta v^\perp$,
where $\delta q_i^0,\, \delta v_i^0\parallel l_0$, and
$\delta q_i^\perp,\, \delta v_i^\perp\perp l_0$ for $i=1,\dots,N$.
Then it is true that

$$
\max\left\{\frac{||\delta q^0||}{||\delta q||},\, \frac{||\delta v^0||}
{||\delta v||}\right\}<\delta_0=\delta_0(\epsilon_0)<<1.
\tag 8.3
$$
Here the small number $\delta_0=\delta_0(\epsilon_0)$ depends on 
$\epsilon_0$ in such a way that it can be made arbitrarily small by selecting 
$\epsilon_0$ small enough.

\endproclaim

\medskip

\subheading{\bf Proof} We observe first that --- since both $J$ and $\Cal S$
are locally flow-invariant --- for any normal vector
$n=(z,w)\in\Cal T_y\bold M$
of $J$ (of $\Cal S$) it is automatically true that $z\perp v$, $w\perp v$,
see Remark 7.9, particularly (7.10). (We always use the notation $y=(q,v)$.)
We note that the orthogonality $w\perp v$ is automatic, for any velocity 
variation $w$ of $v$ is necessarily perpendicular to $v$, due to the energy
normalization $||v||=1$ in the phase space. The reason why the properties
(1)---(3) above can, indeed, be achieved for a smooth foliation
$J=\bigcup_{i\in I}F_i$ ($\text{dim}F_i=1$) is as follows: The unit tangent
vectors $\tau=(\delta q,\, \delta v)\in\Cal T_yJ$ of the curves $F_i$ (yet to
be constructed) have to be, first of all, perpendicular to the normal vector
$n=(z,w)\in\Cal T_y\bold M$ of $J$ at $y=(q,v)\in J$, the vectors $\delta q$,
$\delta v$ have to come from the ortho-complement space $v^\perp$ and, at the
same time, the angles subtended by the vectors $\tau$ and the subspaces
$\Cal T_y\Cal S$ ($y\in J$) have to be separated from zero (uniform
transversality). These things can be achieved simultaneously, according to the
inequalities (7.8) and (7.11). The quadratic forms 
$Q(\tau)=\langle\delta q,\delta v\rangle$ are indefinite on the space
$v^\perp\times v^\perp$ and, consequently, the positivity condition in (2)
still allows a non-empty, open region in 
$\Cal T_yJ\bigcap\left(v^\perp\times v^\perp\right)$ for the unit tangent
vector $\tau=(\delta q,\, \delta v)\in\Cal T_yF_i$. The uniform transversality
of (1) is automatically achieved by the fact that 
$Q(\tau)/\Vert\delta q\Vert^2$ is separated from
zero in (2). The last requirement (8.3) is independent of the former ones,
and it still leaves a non-empty, open set of unit vectors $\tau$ for the
construction of the leaves $F_i$. By integrating the arising, smooth
distribution $\tau(y)$ ($||\tau(y)||=1$, $y\in G\subset J$, 
$G\ne\emptyset$ is an open subset of $J$) on a small, open subset $G$ of
$J$, we obtain a smooth foliation $G=\bigcup_{i\in I}F_i$. Finally, the
original separating manifold is to be replaced by $G$. \qed

\medskip

\subheading{\bf The expansion rate}

\medskip

Consider a non-zero tangent vector
$\tau(0)=\left(\delta q(0),\, \delta v(0)\right)\in\Cal T_yF_i$ of the leaf
$F_i$ at $y\in F_i$. Let us focus on the time-evolution of the so called
infinitesimal Lyapunov function 
$Q(t)=Q(\tau(t))=\langle\delta q(t),\, \delta v(t)\rangle$
($t\ge0$, $\tau(t)=(DS^t)(\tau(0))$) along the non-singular forward orbit
$S^{(0,\infty)}y$, $y\in F_i\setminus\Cal S_J$. The time-evolution of
$\tau(t)=\left(\delta q(t),\, \delta v(t)\right)$ is governed by the equations
(7.1) and (7.3). From those equations we immediately derive the following
time-evolution equations for $Q(t)$ along $S^{(0,\infty)}y$:

$$
\frac{d}{dt}Q(t)=||\delta v(t)||^2\quad\text{ (between collisions)},
\tag 8.4
$$

$$
\aligned
Q(t+0)-Q(t-0)&=2\cos\phi\left\langle RV^*KV\delta q(t-0),\,
R\delta q(t-0)\right\rangle \\
&=2\cos\phi\left\langle KV\delta q(t-0),\,
V\delta q(t-0)\right\rangle\ge 0
\endaligned
\tag 8.5
$$
if a collision takes place at time $t$. In (8.5) we used the well known fact
that $K\ge0$, i. e. the semi-dispersing property. The first consequence of
(8.4)--(8.5) is that the infinitesimal Lyapunov function $Q(t)$ is
non-decreasing in $t$. By the first equation of (7.3), the function
$\left\Vert\delta q(t)\right\Vert^2$ is continuous in $t$ even at
collisions. Its time-derivative between collisions is obtained from the
second equation of (7.1):

$$
\frac{d}{dt}\left\Vert\delta q(t)\right\Vert^2=2\langle\delta q(t),\,
\delta v(t)\rangle=2Q(t).
\tag 8.6
$$
We note that, according to the canonical identification of the tangent
vectors of $\bold Q$ along any trajectory (see \S2 of [K-S-Sz(1990)],
more precisely, the paragraph of that section beginning at the bottom of
p. 538 and ending at the top of p. 539), in the second equation of (7.3)
any tangent vector $w\in\Cal T_{x^-}\bold Q$ gets identified with
$Rw\in\Cal T_{x^+}\bold Q$ ($x^-=S^{t-0}y$, $x^+=S^{t+0}y$), and after
this customary and natural identification the second line of (7.3) turns
into

$$
\delta v^+=\delta v^-+2\cos\phi V^*KV\delta q^-.
\tag 8.7
$$
We recall that the symmetric operator $V^*KV$ in (8.7) is nonnegative. The key
to the understanding of the rate of increase of the function
$\left\Vert\delta q(t)\right\Vert^2$ is that the initial velocity variation
vector $\delta v(0)$ (a component of 
$\tau(0)=\left(\delta q(0),\, \delta v(0)\right)\in\Cal T_yF_i$) can be
obtained as $\delta v(0)=B(0)\delta q(0)$ in such a way that the positive,
symmetric operator $B(0):\; v(0)^\perp\to v(0)^\perp$ is the second
fundamental form of a strictly convex, local orthogonal manifold
$\Sigma\ni y$, and

$$
B(0)\ge c_0 I,
\tag 8.8
$$
see (2) in Proposition 8.2. Denote by $B(t)$ the positive definite second
fundamental form of $S^t\Sigma$ at the point $S^ty$, $t\ge0$. The 
time-evolution of the operators $B(t)$ is governed by the equations (i)--(ii)
of Proposition 2.3 in [K-S-Sz(1990)], see also the last displayed formula
of \S1 in [S-Ch(1982)], or formula (2) in \S3 of [S-Ch(1987)]:

$$
B(t+s)^{-1}=B(t)^{-1}+s\cdot I
\tag 8.9
$$
for $t,\, s\ge0$, provided that $S^{[t,t+s]}y$ is collision free, and

$$
RB(t+0)R=B(t-0)+2\cos\phi V^*KV
\tag 8.10
$$
for a collision at time $t$. From $\delta v(0)=B(0)\delta q(0)$ we obtain

$$
\frac{d}{dt}\delta q(t)=\delta v(t)=B(t)\delta q(t),
\tag 8.11
$$
thus

$$
\delta q(t)=\delta q(0)+\int_0^t B(s)\delta q(s)ds
\tag 8.12
$$
for all $t\ge0$. The equations (8.8)--(8.10) and $V^*KV\ge0$ imply that

$$
B(t)\ge\frac{c_0}{1+c_0t}I\text{ for all } t\ge0.
\tag 8.13
$$
Therefore, 

$$
\aligned
Q(t)=\langle\delta q(t),\, \delta v(t)\rangle
&=\langle\delta q(t),\, B(t)\delta q(t)\rangle \\
\ge\left\langle\delta q(t),\, \frac{c_0}{1+c_0t}\delta q(t)\right\rangle
&=\frac{c_0}{1+c_0t}\Vert\delta q(t)\Vert^2,
\endaligned
$$
so by (8.6) we have

$$
\frac{d}{dt}\Vert\delta q(t)\Vert^2\ge\frac{2c_0}{1+c_0t}
\Vert\delta q(t)\Vert^2,
\tag 8.14
$$
that is,

$$
\frac{d}{dt}\log \Vert\delta q(t)\Vert^2\ge\frac{2c_0}{1+c_0t}.
\tag 8.15
$$
By integration we immediately obtain

$$
\frac{\Vert\delta q(t)\Vert}{\Vert\delta q(0)\Vert}\ge 1+c_0t.
\tag 8.16
$$

\medskip

\subheading{\bf Remark 8.17} It might be interesting to contemplate a bit
about the fact that the lower estimation for $||\delta q(t)||$ is only linear
in $t$. Apparently, the reason is that along a considered forward trajectory
$S^{[0,\infty)}y\subset \bar U_{\epsilon_0}\left(L(l_0)\right)$ the free
path length is actually unbounded, and this fact is known to have the
potential for spoiling any better estimation.

\medskip

For us the utmost important inequality is the lower estimation (8.16)
for the growth of $||\delta q(t)||$. The only shortcoming of (8.16) is
that in the following proof we will need a sufficiently large coefficient
of $t$ on the right-hand-side, instead of just $c_0$. However, this goal can
be achieved as the proof of the following corollary shows.

\medskip

\proclaim{\bf Corollary 8.18} For an arbitrarily big constant $c_1>>1$
one can find a non-empty, open subset $G\subset J$ (and can rename $G$ as
$J$ afterward, as we always do) with the property that the foliation
$G=\bigcup_{i\in I}F_i$ of $G$ (given by the constructive proof of 
Proposition 8.2) can actually be constructed in such a way that the
dilation constant $c_0$ in (2) (and in (8.16))
is replaced by the given number $c_1$.
\endproclaim

\medskip

\subheading{\bf Proof}
Select and fix a phase point $y_0\in J$ with a non-singular
forward orbit $S^{(0,\infty)}y_0$. First appropriately construct the unit
tangent vector $\tau=(\delta q,\, \delta v)\in\Cal T_{y_0}F_i$ of the curve
$F_i$ (to be constructed). The constant $c_0$ in (2) is determined by the
local geometry of $J$ around $y_0$, so it can be chosen to be the same for all
$y\in G_1$ in a suitable neighborhood $G_1$ of $y_0$ in $J$. Now select a
unit tangent vector 
$\tau=(\delta q,\,\delta v)=(\delta q(0),\,\delta v(0))\in\Cal T_{y_0}\bold M$
by using the constructive proof of Proposition 8.2, and also select a time
moment $t_0>c_1/c_0^2$ so that $S^{t_0}y_0\in\partial\bold M$,
i. e. $t_0$ is a moment of collision on the forward orbit
$S^{(0,\infty)}y_0$. Choose a very small $\epsilon_0'>0$ so that
$S^{(t_0,t_0+\epsilon_0']}y_0\cap\partial\bold M=\emptyset$. By (8.16)

$$
\frac{\Vert\delta q(t_0+\epsilon_0')\Vert}{\Vert\delta q(0)\Vert}>
c_0t_0.
\tag 8.19
$$
(Here, as always, we use the notation 
$\tau(t)=(\delta q(t),\, \delta v(t))=(DS^t)(\tau(0))$.) Clearly, there is an
absolute constant $c_2>0$ such that the inequality

$$
\frac{\left|\left\langle\delta q_i(t_0-0)-\delta q_j(t_0-0),\,
l_0^\perp\right\rangle\right|}{\Vert\delta q(t_0-0)\Vert}>c_2
\tag 8.20
$$
can be achieved by suitably selecting the initial (unit) tangent vector
$\tau(0)\in\Cal T_{y_0}\bold M$. Here $i$ and $j$ are the labels of the two
disks colliding at time $t_0$ on $S^{(0,\infty)}y_0$. The reason why (8.20) 
can be achieved is that this inequality defines a non-empty, open cone in terms
of $\delta q(t_0-0)$, and the mapping $\delta q(0)\mapsto\delta q(t_0-0)$ is a
linear bijection between $v(0)^\perp$ and $v(t_0-0)^\perp$ for any given family

$$
\left\{(\delta q(0),\, B\delta q(0))\big|\; \delta q(0)\perp v(0)\right\}
$$
of tangent vectors, where $B\ge0$ and

$$
\left(\delta q(t_0-0),\,\delta v(t_0-0)\right)=DS^{t_0-0}\left
(\delta q(0),\,B\delta q(0)\right).
$$
A consequence of (8.20) is that we obtain the estimation

$$
B(t_0+\epsilon_0')\ge c_3\cdot I
\tag 8.21
$$
of type (8.8) with an absolute constant $c_3>0$. We can assume that the
original $c_0$ is smaller than $c_3$. Then the whole proof of Proposition 8.2
can be repeated for the sub-manifold $S^{t_0+\epsilon_0'}(G_2)$ with some
small, open neighborhood $G_2$ of $y_0$ in $G_1$
($y_0\in G_2\subset G_1\subset J$).
The arising foliation $G_2=\bigcup_{i\in I}F_i$ will enjoy the property
that the $||\delta q||$-expansion rate between $t=0$ and $t=t_0+\epsilon_0'$
is greater than $c_0t_0$ (see also (8.19)), while this rate between
$t_0+\epsilon_0'$ and $t$ ($t>>t_0$) is at least
$c_3(t-t_0-\epsilon_0')\approx c_3t>c_0t$. However, the product of these two
lower estimations of the $||\delta q||$-expansion rates is equal to
$c_0^2t_0t$, which is greater than $c_1t$ by the selection of $t_0$
($t_0>c_1/c_0^2$). This concludes the proof of the corollary. \qed

\medskip

\subheading{The invariant cone field}

\medskip

Now let us pay attention to the cones defined by the inequality (8.3)
and the convexity condition $\langle\delta q,\, \delta v\rangle>0$.
For such tangent vectors $\tau=(\delta q,\, \delta v)$ use the usual
decomposition $\delta q=\delta q^0+\delta q^\perp$, 
$\delta v=\delta v^0+\delta v^\perp$, just as in (3) of Proposition 8.2.

Along a forward orbit 
$S^{[0,\infty)}y\subset\bar U_{\epsilon_0}\left(L(l_0)\right)$ the dilation
effect of the billiard flow {\it between two consecutive collisions} is
dramatically different for the tangent vectors $\tau=(\delta q,\delta v)$
with $\delta q^\perp=\delta v^\perp=0$ (but still 
$\langle\delta q,\, \delta v\rangle>0$, as always in our considerations)
and for the tangent vectors $\tau=(\delta q,\delta v)$ with
$\delta q^0=\delta v^0=0$. By this dramatic difference we mean the following
fact: Let $y\in J$,
$S^{[0,\infty)}y\subset\bar U_{\epsilon_0}\left(L(l_0)\right)$,
$0<t_1<t_2$ any two time moments for which $S^{t_1}y\not\in\partial\bold M$,
$S^{t_2}y\not\in\partial\bold M$, and the non-singular orbit segment
$S^{[t_1,t_2]}y$ has a connected collision graph. Assume that 
$\tau(t_1)=\left(\delta q(t_1),\, \delta v(t_1)
\right)\in\Cal T_{S^{t_1}y}\bold M$,
$\rho(t_1)=\left(\delta\tilde q(t_1),\, \delta\tilde v(t_1)
\right)\in\Cal T_{S^{t_1}y}\bold M$ are two tangent vectors of 
$\bold M$ at $S^{t_1}y$ with the usual convexity property
$Q(\tau(t_1))>0$, $Q(\rho(t_1))>0$. Assume, finally, that 

$$
\aligned
\tau^0(t_1)&=:\left(\delta q^0(t_1),\, \delta v^0(t_1)\right)=(0,0), \\
\rho^\perp(t_1)&=:\left(\delta\tilde q^\perp(t_1),\,
\delta\tilde v^\perp(t_1)\right)=(0,0).
\endaligned
$$
There is a constant $\Lambda>1$ (independent of $y\in J$, $t_1$, $t_2$, 
$\tau(t_1)$, and $\rho(t_1)$, depending only on $N,\,m_1,\dots,m_N$, and
$\epsilon_0$) such that

$$
\frac{||\tau(t_2)||}{||\tau(t_1)||}\div
\frac{||\rho(t_2)||}{||\rho(t_1)||}\ge\Lambda.
\tag 8.22
$$
The reasons why (8.22) holds true are as follows:

\medskip

(1) All collision normal vectors of the trajectory segment $S^{[t_1,t_2]}y$
are almost parallel or orthogonal to the fixed lattice vector $l_0$.
(The angular deviation from the exact parallelity or orthogonality is less
than $\epsilon_0$.) This means that the components $\delta q$, $\delta v$
of the tangent vectors $DS^{t-t_1}\left(\rho(t_1)\right)$ which are almost
parallel to $l_0$ will again be taken into such vectors by the orthogonal
reflection part $R(\,.\,)$ (see (7.3)) of the linearization of the flow at any
collision $S^ty$ ($t_1<t<t_2$), and an analogous statement holds true for the
components $\delta q$, $\delta v$ of the tangent vectors 
$DS^{t-t_1}\left(\tau(t_1)\right)$ which are almost perpendicular to $l_0$.

\medskip

(2) The ``scattering effect'' of the linearized billiard flow at a collision 
$S^ty$ ($t_1<t<t_2$) (i. e. the term $2\cos\phi RV^*KV\delta q^-$ in (7.3))
is almost perpendicular to $l_0$, and this vector is of higher order of 
magnitude for the $\tau$ vectors than for the $\rho$ vectors. Actually, the
ratio of these two effects tends to infinity as $\epsilon_0\to 0$.

\medskip

A direct consequence of the above arguments is

\medskip

\proclaim{\bf Proposition 8.23} Use all of the above assumptions and notations,
that is, that the collision graph of the non-singular orbit
segment $S^{[t_1,t_2]}y$ is connected, $y\in J$, $0<t_1<t_2$,
$S^{t_1}y\not\in\partial\bold M$, $S^{t_2}y\not\in\partial\bold M$.
We claim that the cone field $\Cal C(z)$ ($z=S^{t_1}y$, $y\in J$) defined by 
(8.3) and the convexity condition is invariant under the linearization of the 
billiard map $DS^{t_2-t_1}$, that is, for any
tangent vector $\tau(t_1)\in\Cal T_{S^{t_1}y}\bold M$ ($y\in J$) with
$Q\left(\tau(t_1)\right)>0$ and (8.3) it is true that
$Q\left(\tau(t_2)\right)>0$ and (8.3) still holds for $\tau(t_2)$.

\endproclaim

\medskip

Finally, let us investigate the extent to which the inequality (8.3) can be
spoiled by the free flight between collisions. Use all the notations from
above. Consider a tangent vector
$\tau(t_1+0)=\left(\delta q(t_1+0),\,
\delta v(t_1+0)\right)\in\Cal C\left(S^{t_1}y\right)$ 
($y\in J\setminus\Cal S_J$, $t_1>0$ is a moment of collision on
$S^{(0,\infty)}y$, $S^{(0,\infty)}y\subset\bar
U_{\epsilon_0}\left(L(l_0)\right)$) of the cone field $\Cal C$. Let,
furthermore, $t$ be a positive number with $t<t_2-t_1$ ($0<t_1<t_2$
are the time moments of two consecutive collisions on $S^{(0,\infty)}y$).
Then we claim

\medskip

\proclaim{\bf Proposition 8.24} Use all the above notations. The inequalities

$$
\frac{\Vert\delta q^0(t_1+t)\Vert}{\Vert\delta q(t_1+t)\Vert}<
\sqrt{2}\delta_0,
\tag 8.25
$$

$$
\frac{\Vert\delta v^0(t_1+t)\Vert}{\Vert\delta v(t_1+t)\Vert}<\delta_0
\tag 8.26
$$
hold true.
\endproclaim

\medskip

\subheading{\bf Proof} By the time-evolution equations (7.1) (which hold
separately for $\delta q^0$, $\delta v^0$ on one hand, and 
$\delta q^\perp$, $\delta v^\perp$ on the other hand) we have that
$\delta v(t_1+t)=\delta v(t_1+0)$, 
$\delta v^0(t_1+t)=\delta v^0(t_1+0)$, thus (8.26) is obviously true.
The convexity condition

$$
\left\langle\delta q(t_1+0),\, \delta v(t_1+0)\right\rangle>0
$$
immediately provides the inequality

$$
\left\Vert\delta q(t_1+0)+t\delta v(t_1+0)\right\Vert>
\sqrt{\left\Vert\delta q(t_1+0)\right\Vert^2+t^2\cdot
\left\Vert\delta v(t_1+0)\right\Vert^2}.
\tag 8.27
$$
By the triangle inequality and by the assumption
$\tau(t_1+0)\in\Cal C\left(S^{t_1}y\right)$ we have that

$$
\left\Vert\delta q^0(t_1+0)+t\delta v^0(t_1+0)\right\Vert<\delta_0\cdot
\left(\left\Vert\delta q(t_1+0)\right\Vert+
t\left\Vert\delta v(t_1+0)\right\Vert\right).
\tag 8.28
$$
Combining the inequalities (8.27)--(8.28) with the trivial inequality

$$
a+b\le\sqrt{2(a^2+b^2)}\quad a,\, b\ge 0,
$$
one gets

$$
\frac{\Vert\delta q^0(t_1+t)\Vert}{\Vert\delta q(t_1+t)\Vert}<
\frac{\delta_0\left(\Vert\delta q(t_1+0)\Vert+
t\Vert\delta v(t_1+0)\Vert\right)}
{\sqrt{\left\Vert\delta q(t_1+0)\right\Vert^2+t^2\cdot
\left\Vert\delta v(t_1+0)\right\Vert^2}}\le\sqrt{2}\delta_0,
$$
which finishes the proof of the proposition. \qed

\medskip

\proclaim{\bf Corollary 8.29. (Corollary of (8.16) and propositions 
8.23--8.24)} For any tangent vector
$\tau=(\delta q(0),\, \delta v(0))\in\Cal C(y)$
($y\in J\setminus\Cal S_J$, $S^{(0,\infty)}y\subset\bar
U_{\epsilon_0}\left(L(l_0)\right)$) it is true that

$$
\lim_{t\to\infty}\frac{\Vert\delta q^0(t)\Vert}{\Vert\delta q(t)\Vert}=0.
$$
\endproclaim

\medskip

\subheading{\bf Frequency of collisions (Frequency of singularities)}

\medskip

Denote by $\#\left(S^{[a,b]}y\right)$ the number of collisions on the
non-singular trajectory segment $S^{[a,b]}y$. Assume that the non-degeneracy 
condition of Corollary 1.1 of [B-F-K(1998)] holds true at all phase points
$x\in\partial\bold M$ lying close enough to the limiting set $L(l_0)$.
This condition at a boundary phase point $x\in\partial\bold M$ essentially
means that the spatial angle subtended by $\text{int}\bold M$ at $x$
is positive. It is easy to see that this positive-angle condition can only be
violated if either

(i) $2r\vert H_i\vert$ is equal to the length $\Vert l_0\Vert$ of the closed
geodesic of $\Bbb T^2$ in the direction of the vector $l_0$, (The vector
$l_0$ is supposed to be non-divisible in $\Bbb Z^2$.);

\noindent
or

(ii) $2rk$ is equal to the width of the torus $\Bbb T^2$ in the direction
of the perpendicular vector $l_0^\perp$.

Recall that $k$ denotes the number of different groups of disks $H_i$, see 
the paragraph right before (6.2). The width of $\Bbb T^2$ in the direction
of $l_0^\perp$ is, by definition, equal to the length of the shortest
vector in the orthogonal projection of $\Bbb Z^2$ onto the line spanned by
$l_0^\perp$. (This length is just the reciprocal of $\Vert l_0\Vert$.)

There are only countably many values of the radius $r$ for which either
(i) or (ii) is true, and those exceptional values may be discarded without
narrowing the scope of our Theorem.

It follows from Corollary 1.1 of [B-F-K(1998)] that there exists a constant
$c_4'=c_4'(N,r,m_1,\dots,m_N)>0$ (depending only on the geometry of the hard
disk system) such that $\#\left(S^{[a,a+1]}y\right)\le c_4'$ for all
non-singular orbit segments $S^{[a,a+1]}y$. Consequently, there exists another
constant $c_4=c_4(N,r,m_1,\dots,m_N)>0$ such that

$$
\#\left(S^{[a,a+t]}y\right)\le c_4\text{max}\{t,\,1\}
\tag 8.30
$$
for all non-singular trajectory segments $S^{[a,a+t]}y$.

\medskip

\subheading{Sinai's idea: ``Expansion prevails over chopping''
(Finishing the proof of the Theorem)}

\medskip

Take a large constant $c_1>>1$ and, by using Corollary 8.18 above,
consider a smooth foliation $J=\bigcup_{i\in I}F_i$ by curves $F_i$
fulfilling all conditions listed in Proposition 8.2 in such a way that
the expansion constant $c_0$ in (8.16) is actually the large constant $c_1$.
Later in the proof we will see how large the constant $c_1$ should actually
be chosen in order that the whole proof of the Theorem works. Pick up a single
curve $F_{i_0}=F_0\subset J$ of the foliation $J=\bigcup_{i\in I}F_i$.
On the curve $F_0$ itself and on the connected components of its forward
images $S^t(F_0)$ we will be measuring the distances by using the so called
$z$-distance introduced by Chernov and Sinai (cf. Lemma 2 and the preceding
paragraph in \S4 of [S-Ch(1987)]) defined as follows:

$$
z(y_1,y_2)=:\int_{y_1}^{y_2}\Vert dq\Vert
\tag 8.31
$$
for points $y_1,\, y_2$ of a connected component $\gamma$ of the image
$S^t(F_0)$. The integral in (8.31) is taken on the segment of $\gamma$
connecting $y_1$ and $y_2$. Set

$$
S(t)=\big\{y\in F_0\big|\; S^{[0,t]}y
\text{ contains at least one singular collision}\big\}.
\tag 8.32
$$
By (1) of Proposition 8.2 we see that $F_0$ intersects any future singularity
manifold $\Cal S$ in at most one point, and the number of such singularity
manifolds until time $t$ is at most $c_4t$ by (8.30), so we get the following
upper estimation for the cardinality $k(t)$ of the set $S(t)$:

$$
k(t)=:\left|S(t)\right|\le c_4t\text{ for all }t\ge 1.
\tag 8.33
$$
(We are only interested in large values of $t$.) Let

$$
F_0\setminus S(t)=\cup_{p=1}^{k(t)+1}I_p^{(t)}
$$
be the decomposition of the open set $F_0\setminus S(t)$ into its connected
components. Select a positive constant $c_5$ (for its actual value, see
below), and define

$$
B(t)=\bigcup\left\{I_p^{(t)}\Big|\; \left|I_p^{(t)}\right|_z<c_5/t \right\}.
\tag 8.34
$$
Here the length $\left|I_p^{(t)}\right|_z$ of $I_p^{(t)}$ is measured by using
the $z$-metric of (8.31). From (8.33)--(8.34) we obtain the estimation

$$
\mu_z\left(B(t)\right)<\frac{c_5}{t}\cdot c_4t=c_4c_5<
\frac{1}{2}\mu_z(F_0),
\tag 8.35
$$
as long as the constant $c_5>0$ is selected so that $c_5<\mu_z(F_0)/(2c_4)$. 
Here $\mu_z$ is the Lebesgue measure on the curve $F_0$ defined by the 
distance parametrization $z$ from (8.31). We recall that the foliation 
$J=\bigcup_{i\in I}F_i$ (and, consequently, the chosen curve $F_0$, as well)
depends on the constant $c_1$. For any component
$I_p^{(t)}\subset G(t)=:F_0\setminus\overline{B(t)}$ we have
$\left|I_p^{(t)}\right|_z\ge c_5/t$, and by (8.16) (with $c_0$ replaced
by $c_1$) we get the estimation

$$
\mu_z\left(S^t\left(I_p^{(t)}\right)\right)>c_1t\cdot\frac{c_5}{t}
=c_1c_5.
\tag 8.36
$$
Use the shorthand notation $\gamma_p=S^t\left(I_p^{(t)}\right)$ for any
$I_p^{(t)}$, $I_p^{(t)}\subset G(t)$. By the invariance of the cone field
$\Cal C(z)$ along any trajectory 
$S^{(0,\infty)}y\subset\bar U_{\epsilon_0}\left(L(l_0)\right)$ (with the
additional features $\lim_{t\to\infty}d(S^ty,\, L(l_0))=0$,
$y\in J\setminus\Cal S_J$), see particularly Corollary 8.29, it is true that
the integral $\int_{\gamma_p}||dq||$ is asymptotically the same as
$\int_{\gamma_p}\left|\langle dq,\, l_0^\perp\rangle\right|$ and, accordingly,
$\int_{\gamma_p}||dv||$ is also asymptotically the same as
$\int_{\gamma_p}\left|\langle dv,\, l_0^\perp\rangle\right|$, where 
$l_0^\perp\in\Bbb R^2$ is a formerly selected unit vector perpendicular to the
lattice vector $l_0$ defining $L(l_0)$. What is even more, the scattering
property of the hard disk system along the studied orbits 
$S^{(0,\infty)}y$ ($y\in F_0$) is such that there exists a constant
$c_6=c_6(N,r,m_1,\dots,m_N)$ (again depending only on the geometry of the
hard disk system) such that

$$
\int_{\gamma_p}\left|\langle dv,\, l_0^\perp\rangle\right|\ge
c_6\cdot\int_{\gamma_p}\Vert dq\Vert>c_1c_5c_6=:100c_7.
\tag 8.37
$$
(In the second inequality we used (8.36).) Use the shorthand 
$c_7=:\dfrac{c_1c_5c_6}{100}$ in (8.37). By reversing the simple dilation
argument based upon (8.16) (with $c_0$ replaced by $c_1$) and leading to
(8.36), we get that for any pair of points
$y_1,\, y_2\in\gamma_p\cap U_{c_7}\left(L(l_0)\right)$
($U_{c_7}\left(L(l_0)\right)$ denotes the open $c_7$-neighborhood
of the compact set $L(l_0)$) it is true that

$$
\int_{y_1}^{y_2}\left|\langle dv,\, l_0^\perp\rangle\right|\le 2c_7
$$
and, consequently,

$$
z\left(S^{-t}y_1,\, S^{-t}y_2\right)<\frac{c_5}{50t}.
\tag 8.38
$$
Set $C(t)=:F_0\cap S^{-t}\left(U_{c_7}\left(L(l_0)\right)\right)$.
An immediate consequence of (8.38) is that

$$
\frac{\mu_z\left(C(t)\cap I_p^{(t)}\right)}
{\mu_z\left(I_p^{(t)}\right)}\le\frac{1}{50}
\tag 8.39
$$
for all $I_p^{(t)}$, $I_p^{(t)}\subset G(t)=:F_0\setminus\overline{B(t)}$.
Consequently,

$$
\frac{\mu_z\left(C(t)\cap G(t)\right)}{\mu_z\left(G(t)\right)}
\le\frac{1}{50},
$$
thus

$$
\mu_z\left(C(t)\right)<\left(\frac{1}{2}+\frac{1}{50}\right)\mu_z(F_0),
\tag 8.40
$$
by also using (8.35). By Proposition 6.5, however, we have that

$$
\lim_{t\to\infty}\mu_z\left(C(t)\right)=\mu_z(F_0),
$$
in contradiction with (8.40). This step completes the proof of the 
non-existence of any separating manifold $J$, thereby finishing the whole
proof of the Theorem. \qed

\bigskip \bigskip

\heading
\S 9. Concluding Remarks
\endheading

\bigskip \bigskip

\subheading{9.1 Irrational Mass Ratio}

\medskip

Due to the natural reduction $\sum_{i=1}^N m_iv_i=0$ (which we always assume),
in \S 1 we had to factorize out the configuration space with respect to
spatial translations: $(q_1,\dots,q_N)\sim(q_1+a,\dots,q_N+a)$ for all
$a\in\Bbb T^2$. It is a remarkable fact, however, that (despite the 
reduction $\sum_{i=1}^N m_iv_i=0$) even without this translation factorization
the system still retains the Bernoulli mixing property, provided that the
masses $m_1,\dots,m_N$ are rationally independent. (We note that dropping the
above mentioned configuration factorization obviously introduces $2$ zero
Lyapunov exponents.) For the case $N=2$ (i. e. two disks)
this was proved in [S-W(1989)] by successfully applying D. Rudolph's theorem
on the B-property of isometric group extensions of Bernoulli shifts [R(1978)].

Suppose that we are given a dynamical system $(M,T,\mu)$ with a probability
measure $\mu$ and an automorphism $T$. Assume that a compact metric group
$G$ is also given with the normalized Haar measure $\lambda$ and left invariant
metric $\rho$. Finally, let $\varphi\colon\; M\to G$ be a measurable map.
Consider the skew product dynamical system $(M\times G,S,\mu\times\lambda)$ 
with $S(x,g)=\left(Tx,\varphi(x)\cdot g\right)$, 
$x\in M$, $g\in G$. We call the system $(M\times G,S,\mu\times\lambda)$ an 
isometric group extension of the base (or factor) $(M,T,\mu)$. (The phrase
``isometric'' comes from the fact that the left translations 
$\varphi(x)\cdot g$ are isometries of the group $G$.) Rudolph's mentioned
theorem claims that the isometric group extension 
$(M\times G,S,\mu\times\lambda)$ enjoys the B-mixing property as long as it is
at least weakly mixing and the factor system $(M,T,\mu)$ is a B-mixing system.

But how do we apply this theorem to show that the system of $N$ hard disks
on $\Bbb T^2$ with $\sum_{i=1}^N m_iv_i=0$ is a Bernoulli flow, even if we
do not make the factorization (of the configuration space) with respect to
spatial translations? It is simple. The base system $(M,T,\mu)$
of the isometric group extension $(M\times G,S,\mu\times\lambda)$ will be the
time-one map of the factorized (with respect to spatial translations) hard 
disk system. The group $G$ will be just the container torus $\Bbb T^2$
with its standard Euclidean metric $\rho$ and normalized Haar measure
$\lambda$. The second component $g$ of a phase point
$y=(x,g)\in M\times G$ will be just the position of the center of the (say) 
first disk in $\Bbb T^2$. Finally, the governing translation
$\varphi(x)\in\Bbb T^2$ is quite naturally the total displacement

$$
\int\Sb 0\endSb\Sp 1\endSp v_1(x_t)dt \qquad (\text{mod }\Bbb Z^2)
$$
of the first particle while unity of time elapses. In the previous sections
the B-mixing property of the factor map $(M,T,\mu)$ has been proved
successfully for typical geometric parameters $(m_1,\dots,m_N;r)$.
Then the key step in proving the B-property of the isometric
group extension $(M\times G,S,\mu\times\lambda)$ is to show that the latter
system is weakly mixing. This is just the essential contents of the
paper [S-W(1989)], and it takes advantage of the assumption of rational
independence of the masses. Here we are only presenting to the reader
the outline of that proof. As a matter of fact, we not only proved
the weak mixing property of the extension $(M\times G,S,\mu\times\lambda)$,
but we showed that this system has in fact the K-mixing property by proving
that the Pinsker partition $\pi$ of $(M\times G,S,\mu\times\lambda)$ is
trivial. (The Pinsker partition is, by definition, the finest measurable
partition of the dynamical system with respect to which the factor system
has zero metric entropy. A dynamical system is K-mixing if and only if
its Pinsker partition is trivial, i. e. it consists of only sets with
measure zero and one, see [K-S-F(1980)].) In order to show that the Pinsker
partition is trivial, in [S-W(1989)] we constructed a pair of measurable
partitions $(\xi^s,\,\xi^u)$ for $(M\times G,S,\mu\times\lambda)$ made up
by open and connected sub-manifolds of the local stable and unstable manifolds,
respectively. It followed by standard methods (see [Sin(1968)]) that the 
partition $\pi$ is coarser than each of $\xi^s$ and $\xi^u$. Due to the 
$S$-invariance of $\pi$, we have that $\pi$ is coarser than
$$
\bigwedge\Sb n\in\Bbb Z\endSb S^n\xi^s\wedge
\bigwedge\Sb n\in\Bbb Z\endSb S^n\xi^u.
\tag 9.2
$$
In the final step, by using now the rational independence of the masses,
we showed that the partition in $(9.2)$ is, indeed, trivial.

\bigskip

\subheading{9.3 The role of Proposition 3.1} By taking a look at \S 3, 
we can see that Proposition 3.1 (with its rather involved algebraic proof)
was only used to prove the so-called Chernov-Sinai Ansatz, an important,
necessary condition of the Theorem on Local Ergodicity. It is exactly the
algebraic proof of Proposition 3.1 that necessitates the dropping of a
null set of geometric parameters $(m_1,\dots,m_N;r)$ in such an implicit
way that for any given $(N+1)$-tuple $(m_1,\dots,m_N;r)$ one cannot tell
(based upon the presented methods) if that $(N+1)$-tuple belongs to the
exceptional null set, or not. This is a pity, indeed, since we cannot make
it sure (for any specified $(N+1)$-tuple $(m_1,\dots,m_N;r)$) that the
billiard flow $\flow$ is ergodic. Thus, it would be really pleasant to
find any other way of proving the Ansatz in order to avoid the necessary
dropping of a null set of parameters. Most experts are absolutely convinced
that, in fact, this exceptional null set is actually empty, i. e. $\flow$
is ergodic for every $(N+1)$-tuple $(m_1,\dots,m_N;r)$).

Without Proposition 3.1, the results of \S 4--8 (the non-existence of
the exceptional $J$-manifold) are easily seen to yield the following, relaxed
version of the Chernov-Sinai Ansatz:

\medskip

\proclaim{Proposition 9.4 (Ansatz, relaxed version)} The closed set
$B\subset\Cal S\Cal R^+$ of phase points $x\in\Cal S\Cal R^+$ with 
non-sufficient semi-orbit $S^{(0,\infty)}x$ is of first category in any
$(2d-3)$-dimensional cell $C$ of $\Cal S\Cal R^+$, which is now equivalent
to saying that $B$ has an empty interior in $C$.
\endproclaim

\bigskip \bigskip

\subheading{Acknowledgement} I would like to express my sincere gratitude to
the referee of the paper for making a large number of very useful remarks and
suggestions which significantly raised the quality of this work.

\bye